\documentstyle[12pt]{article}
\def\mymedskip{\vskip\medskipamount}
\def\mymedbreak{\par \ifdim\lastskip<\medskipamount
  \removelastskip \penalty-100 \mymedskip \fi}
\def\myaftermedspace{\par \ifdim\lastskip<\medskipamount
  \removelastskip \penalty55\mymedskip\fi}
\newcommand{\eop}{{\unskip\nobreak\hfil\penalty50
          \hskip2em\hbox{}\nobreak\hfil$\Box$
          \parfillskip=0pt \finalhyphendemerits=0 \par}}
\newenvironment{proof}%
{\mymedbreak{\noindent\bf Proof:\enspace}}{\eop\myaftermedspace}
{\mymedbreak{\noindent\bf Proof of Theorem #1:\enspace}}{\eop\myaftermedspace}
\newenvironment{remark}%
{\mymedbreak\noindent{\bf Remark:}%
\enspace\rm}%
{\myaftermedspace}
\newtheorem{teor}{Theorem}[section]
\newtheorem{defi}[teor]{Definition}
\newtheorem{examp}[teor]{Example}
\newtheorem{lem}[teor]{Lemma}
\newtheorem{cor}[teor]{Corollary}

\newcommand{\beq}{\begin{equation}}
\newcommand{\eeq}{\end{equation}}
\newcommand{\beql}[1]{\begin{equation} \label{#1}}
\newcommand{\eeql}{\end{equation}}
\newcommand{\beqa}{\begin{eqnarray*}}
\newcommand{\eeqa}{\end{eqnarray*}}
\newcommand{\beqal}[1]{\begin{eqnarray} \label{#1}}
\newcommand{\eeqal}{\end{eqnarray}}
\newcommand{\beqan}{\begin{eqnarray}}
\newcommand{\eeqan}{\end{eqnarray}}
\newcommand{\bpf}{\begin{proof}}
\newcommand{\epf}{\end{proof}}

\newcommand{\cE}{{\cal E}}
\newcommand{\cL}{{\cal L}}
\newcommand{\cO}{{\cal O}}
\newcommand{\Obar}{{\bar{\cal O}}}
\newcommand{\cR}{{\cal R}}
\newcommand{\cS}{{\cal S}}

\newcommand{\bB}{{\bf B}}
\newcommand{\bF}{{\bf F}}
\newcommand{\bR}{{\bf R}}
\newcommand{\bS}{{\bf S}}
\newcommand{\bT}{{\bf T}}

\newcommand{\xb}{\bar{x}}
\newcommand{\yb}{\bar{y}}
\newcommand{\zb}{\bar{z}}

\newcommand{\Db}{\bar{\Delta}}

\newcommand{\PG}{{\rm PG}}
\newcommand{\PGL}{{\rm PGL}}
\newcommand{\PgammaL}{{\rm P\Gamma L}}

\newcommand{\Tr}{{\rm Tr}}

\newcommand{\gf}{\phi}

\newcommand{\ga}{\alpha}
\newcommand{\gb}{\beta}
\newcommand{\gc}{\gamma}
\newcommand{\gd}{\delta}
\newcommand{\gre}{\epsilon}
\newcommand{\gl}{\lambda}

\newcommand{\bE}{{\bf E}}
\newcommand{\bG}{{\bf G}}

\newcommand{\rhoh}{{\hat{\rho}}}

\newcommand{\cH}{{\cal H}}
\newcommand{\cA}{{\cal A}}
\newcommand{\coco}{{coherent configuration}}

\newcommand{\bFq}{{\bF_q}}
\newcommand{\bFqp}{{\bF_q\cup\{\infty\}}}
\newcommand{\bFqmin}{{\bF_{q^2}\setminus\bF_q}}
\newcommand{\Gal}{{\rm Gal}(\bE/\bF)}
\newcommand{\bFqt}{\bF_{q^2}}
\begin{document}
\begin{titlepage}
\title{Association schemes from the action of
$\PGL(2,q)$ fixing a nonsingular conic in $\PG(2,q)$}
\date{\today}
\author{%
Henk D. L.\ Hollmann\\Philips Research Laboratories\\
Prof. Holstlaan 4, 5656 AA Eindhoven\\The Netherlands\\email: {\tt
henk.d.l.hollmann@philips.com}\\
\\Qing Xiang\\Department of Mathematical
Sciences\\University of Delaware\\Newark, DE 19716, USA\\email: {\tt
xiang@math.udel.edu}%
}
\maketitle
\begin{abstract}
The group $\PGL(2,q)$ has an embedding into $\PGL(3,q)$ such that it acts as the group fixing a nonsingular conic
in $\PG(2,q)$. This action affords a coherent configuration $\cR(q)$ on the set $\cL(q)$ of non-tangent lines of
the conic. We show that the relations can be described by using the cross-ratio. Our results imply that the
restrictions $\cR_+(q)$ and $\cR_-(q)$ of $\cR(q)$ to the set $\cL_+(q)$ of secant (hyperbolic) lines and to the
set $\cL_-(q)$ of exterior (elliptic) lines, respectively, are both association schemes; moreover, we show that
the elliptic scheme $\cR_-(q)$ is pseudocyclic.

We further show that the coherent configurations $\cR(q^2)$ with $q$ even allow certain fusions. These provide a
4-class fusion of the hyperbolic scheme $\cR_+(q^2)$, and 3-class fusions and 2-class fusions (strongly regular
graphs) of both schemes $\cR_+(q^2)$ and $\cR_-(q^2)$. The fusion results for the hyperbolic case are known, but
our approach here as well as our results in the elliptic case are new.
\end{abstract}
\end{titlepage}

\section{\label{intro}Introduction}

Let $q$ be a prime power. The 2-dimensional projective linear group $\PGL(2,q)$ has an embedding
into $\PGL(3,q)$ such that it acts as the group $G$ fixing a nonsingular conic
$$\cO=\cO_q=\{(\xi, \xi^2, 1)\mid \xi\in \bF_q\}\cup\{(0,1,0)\}$$
in~$\PG(2,q)$ setwise, see e.g. \cite[p.~158]{hirsch1}. Such a conic consists of $q+1$ points forming an {\em
oval\/}, that is, each line  of $\PG(2,q)$ meets $\cO$ in at most two points. Lines meeting the oval in two
points, one point, or no points at all are called {\em secant\/} (or {\em hyperbolic\/}) lines, {\em tangent}
lines, and exterior (or {\em elliptic\/}) lines, respectively. There is precisely one tangent through each point
of an oval; moreover, if $q$ is even, then all tangents pass through a unique point called the {\em nucleus\/} of
the oval, see e.g.\ \cite[p.~157]{hirsch1}.

It turns out that the group $G$ acts generously transitively on both the set $\cL_+$ of hyperbolic lines and the
set $\cL_{-}$ of elliptic lines. Thus we obtain two (symmetric) association schemes, one on $\cL_+$ and the other
on $\cL_{-}$. We will refer to these schemes as the {\em hyperbolic\/} scheme and the {\em elliptic\/} scheme,
respectively.

Our aim in this paper is to investigate these two association schemes {\em simultaneously\/}. Also investigated
here is a particular fusion of these schemes when $q$ is even. In fact, the hyperbolic and elliptic schemes are
contained in the coherent configuration obtained from the action of $G$ on the set $\cL=\cL_+\cup\cL_{-}$ of all
non-tangent lines of the conic $\cO$, and the fusions of the two schemes arise within a certain fusion of this
coherent configuration.

These schemes as well as their fusions are not completely new, but our treatment will be new. For $q$ even, the
elliptic schemes were first introduced in \cite{henkthesis}, as a family of {\em pseudocyclic\/} association
schemes on nonprime-power number of points. The hyperbolic schemes, and the particular fusion discussed here for
$q$ an even square, turn out to be the same as the schemes investigated in \cite{CD}. The fact that the
particular fusion in the hyperbolic case again produces association schemes has been proved by direct computation
in \cite{note}, by geometric arguments in \cite{eehx}, and by using character theory in~\cite{banstud}. The
fusion schemes for $q$ an even square in the elliptic case seem to be new.

The contents of this paper are as follows. In Section~2 we introduce the definitions and notations that are used
in this paper. Then, in Section~3 we introduce the embedding of $\PGL(2,q)$ as the subgroup $G=G(\cO)$ of
$\PGL(3,q)$ fixing the conic $\cO$ in $\PG(2,q)$.

With each non-tangent line we can associate a pair of points, representing its intersection with $\cO$ in the
hyperbolic case, or its intersection with the extension $\cO_{q^2}$ of $\cO$ to a conic in $\PG(2,q^2)$ in the
elliptic case. In Section~4 we show that the orbits of $G$ on pairs of non-tangent lines can be described with
the aid of the cross-ratio of the two pairs of points associated with the lines. These results are then used to
give (new) proofs of the fact that the group action indeed affords association schemes on both $\cL_+$ and
$\cL_{-}$. Moreover, these results establish the connection between the hyperbolic scheme and the scheme
investigated in \cite{CD}.

In Section~5 we develop an expression to determine the orbit to which a given pair of lines belongs in terms of
their homogeneous coordinates.

From Section~6 on we only consider the case where $q$ is even. In Section~6 we derive expressions for the
intersection parameters of the coherent configuration $\cR(q)$ on the non-tangent lines $\cL$ of the conic $\cO$,
so in particular we obtain expressions for the intersection parameters of both the hyperbolic and elliptic
association schemes {\em simultaneously\/}. We also show that the elliptic schemes are pseudocyclic. In
\cite{pseudo} we will prove that the schemes obtained from the elliptic scheme by fusion with the aid of the Frobenius
automorphism of the underlying finite field $\bF_q$ for $q=2^r$ with $r$ prime are also pseudocyclic.

Then in Section~7 we define a particular fusion of the coherent configuration. The results of the previous
section are used to show that this fusion is in fact again a coherent configuration, affording a four-class
scheme on the set of hyperbolic lines and a three-class scheme on the set of elliptic lines. The parameters show
that the restriction of these schemes to one of the classes produces in fact a {\em strongly regular graph\/},
with the same parameters as the Brouwer-Wilbrink graphs (see \cite{bl}) in the hyperbolic case and as the Metz 
graphs (e.g., \cite{bl}) in the elliptic case. This will be discussed in Section~8. In fact, the graphs are {\em
isomorphic\/} to the Brouwer-Wilbrink graphs (in the hyperbolic case) and the Metz graphs (in the elliptic case).
For the hyperbolic case, this was conjectured in \cite{CD} and proved in \cite{eehx}; for the elliptic case, this
was conjectured for $q=4$ in \cite{henkthesis}, and will be proved for general even $q$ in \cite{hh-iso}.

\section{\label{def} Definitions and notation}

\subsection{Coherent configurations}

Let $X$ be a finite set. A {\em coherent configuration\/} is a collection $\cR=\{R_0, \ldots, R_n\}$ of subsets
of $X\times X$ satisfying the following conditions:
\begin{enumerate}
\item $\cR$ is a partition of $X\times X$;
\item there is a subset $\cR_{\rm diag}$ of $\cR$ which is a partition of the diagonal
$\{(x,x) \mid x \in X\}$;

\item for each $R$ in $\cR$, its transpose $R^{\top} = \{(y,x)\mid (x,y)\in R\}$ is
again in $\cR$;
\item there are integers $p_{i j}^k$, for $0\leq i,j,k\leq n$, such
that for all $(x,y)\in R_k$,
\[|\{z\in X \mid (x,z)\in R_i\; {\rm and}\; (z,y)\in R_j\}|=p_{ij}^k.\]
\end{enumerate}
The numbers $p_{ij}^k$ are called the {\em intersection parameters\/} of the coherent
configuration.

Each relation $R_i$ can be represented by its {\em adjacency matrix\/} $A_i$, a matrix
whose rows and columns are both indexed by $X$ and
\[A_i(x,y)=\left\{
        \begin{array}{ll}
                1, & \mbox{if $(x,y)\in R_i$};\\
                0, &\mbox{otherwise}.
        \end{array}
                \right.
\]
In terms of these matrices, and with $I$, $J$ denoting the identity matrix and the all-one matrix, respectively,
the axioms can be expressed in the following form:
\begin{enumerate}
\item
$A_0+A_1+\cdots +A_n=J$;
\item $\displaystyle\sum_{i=0}^m A_i = I$, where $R_{\rm diag} = \{R_0, \ldots, R_m\}$;
\item for each $i$, there exists $i^*$ such that $A_i^{\top}=A_{i^*}$;
\item for each $i,j \in \{0,1, \ldots, n\}$, we have
\[ A_iA_j=\sum_{k=0}^n p_{ij}^k A_k.\]
\end{enumerate}

As a consequence of Properties 2 and 4, the span 
of the matrices $A_0,A_1,\cdots ,A_n$ over the complex numbers is an algebra. It follows from Property 3 that
this algebra is semi-simple, and so is isomorphic to a direct sum of matrix algebras over the complex numbers.

The sets $Y\subseteq X$ such that $\{(y,y)\mid y\in Y\}\subset \cR$ are called the {\em fibres\/} of $\cR$;
according to Property 2, they form a partition of $X$. The coherent configuration is called {\em homogeneous\/}
if there is only one fibre. In that case one usually numbers the relations of $\cR$ such that $R_0$ is the
diagonal relation.

\vskip5pt
\noindent{\bf Remark:} The existence of the numbers $p^k_{d,k}$ and $p^k_{k,d}$ for all diagonal
relations $R_d\in\cR_{\rm diag}$ implies that for each relation $R_k\in\cR$ there are fibres $Y,Z$ such that
$R_k\subseteq Y\times Z$.

\vskip5pt

A coherent configuration is called {\em symmetric\/} if all the relations are symmetric. As a consequence of the
above remark, a symmetric coherent configuration is homogeneous. Usually, a symmetric coherent configuration is
called a {\em (symmetric) association scheme\/}. In this paper, we will call a coherent configuration {\em weakly
symmetric\/} if the restriction of the coherent configuration to each of its fibres is symmetric, that is, each
of its fibres carries an association scheme.

A {\em fusion\/} of a \coco\ $\cR$ on $X$ is a \coco\ $\cS$ on $X$ where each relation $S\in\cS$ is a union of
relations from~$\cR$.

As a typical example of coherent configuration, if $G$ is a permutation group on a finite set $X$, then the
orbits of the induced action of $G$ on $X\times X$ form a coherent configuration; it is homogeneous precisely
when $G$ is {\em transitive\/}, and an association scheme if and only if $G$ acts {\em generously transitively\/}
on $X$, that is, for all $x, y\in X$, there exists $g\in G$ such that $g(x)=y$ and $g(y)=x$. The coherent
configuration is weakly symmetric precisely when $G$ is generously transitive on each of its orbits on $X$.
\subsection{Association schemes}
In the case of an association scheme, Properties 2 and 3 can be replaced by the stronger
properties: \\
2'. $A_0=I$;
and\\
3'. each $A_i$ is symmetric.\\
As a consequence of these properties, the matrices $A_0=I, A_1, \ldots, A_n$ span an algebra $\cA$ over the reals
(which is called the {\em Bose-Mesner algebra} of the scheme). This algebra has a basis $E_0,E_1,\cdots , E_n$
consisting of primitive idempotents, one of which is $\frac{1}{|X|}J$. So we may assume that $E_0=\frac {1}
{|X|}J$. Let $\mu_i={\rm rank}\;E_i$. Then
$$\mu_0=1,\; \mu_0+\mu_1+\cdots +\mu_n=|X|.$$
The numbers $\mu_0,\mu_1,\ldots ,\mu_n$ are called the {\em multiplicities} of the scheme.

Define $P=\left(P_j(i)\right)_{0\le i,j\le n}$ (the {\em first} eigenmatrix) and
$Q=\left(Q_j(i)\right)_{0\le i,j\le n}$
(the {\em second} eigenmatrix) as the $(n+1)\times (n+1)$ matrices with rows and columns indexed by
$0,1,2,\ldots, n$ such that
\[(A_0, A_1, \ldots, A_n)=(E_0, E_1, \ldots, E_n)P,\]
and
\[|X|(E_0, E_1, \ldots, E_n)=(A_0, A_1, \ldots, A_n)Q.\]
Of course, we have
\[P=|X|Q^{-1}, \;\; Q=|X|P^{-1}.\]
Note that $\{P_j(i)\ |\ 0\le i\le n\}$ is the set of eigenvalues of $A_j$ and the zeroth row and column of $P$
and $Q$ are as indicated below.
\[P=\left(\matrix{
1&v_1&\cdots&v_n\cr 1\cr \vdots & & \cr 1 }\right),\;\;
Q=\left(\matrix{ 1&\mu_1&\cdots&\mu_n\cr 1\cr \vdots & & \cr 1 }\right)\]
The numbers $v_0,v_1,\ldots ,v_n$ are called the {\em valencies (or degrees)} of the scheme.

\begin{examp}\label{cycexamp}
We consider {\it cyclotomic schemes} defined as follows. Let $q$ be a prime power and let $q-1=ef$ with $e>1$. Let $C_0$ be the subgroup of the multiplicative group of $\bF_q$ of index $e$, and let $C_0,C_1,\ldots ,C_{e-1}$ be the cosets of $C_0$. We require $-1\in C_0$. Define $R_0=\{(x,x) : x\in \bF_q\}$, and for $i\in \{1,2,\ldots ,e\}$, define $R_i=\{(x,y)\mid x,y\in \bF_q, x-y\in C_{i-1}\}$. Then $(\bF_q, \{R_i\}_{0\leq i\leq e})$ is an $e$-class symmetric association scheme. The intersection parameters of the cyclotomic scheme are related to the cyclotomic numbers (\cite[p.~25]{st}). Namely, for $i,j,k\in \{1,2,\ldots ,e\}$, given $(x,y)\in R_k$,
\begin{equation}\label{cycparam}
p_{ij}^k=|\{z\in \bF_q\mid x-z\in C_{i-1}, y-z\in C_{j-1}\}|=|\{z\in C_{i-k}\mid 1+z\in C_{j-k}\}|.
\end{equation}
The first eigenmatrix $P$ of this scheme is the following $(e+1)$ by $(e+1)$ matrix (with the rows of $P$ arranged in a certain way)
$$P=\left(\matrix{
1&f&\cdots&f\cr
1\cr
\vdots & & P_0\cr
1
}\right)$$
with $P_0=\sum_{i=1}^{e}\eta_iC^i$, where $C$ is the $e$ by $e$ matrix:
$$C=\left(\matrix{
& 1\cr & & 1\cr & & & \ddots\cr & & & & 1\cr 1}\right)$$ and $\eta_i=\sum_{\beta\in C_i}\psi(\beta)$, $1\leq
i\leq e$, for a fixed nontrivial additive character $\psi$ of $\bF_q$.
\end{examp}

Next we introduce the notion of a pseudocyclic association scheme.

\begin{defi}
Let $(X, \{R_i\}_{0\leq i\leq n})$ be an association scheme. We say that $(X, \{R_i\}_{0\leq i\leq n})$ is {\em
pseudocyclic} if there exists an integer $t$ such that $\mu_i=t$ for all $i\in \{1,\cdots , n\}$.
\end{defi}

The following theorem gives combinatorial characterizations for an association scheme to be
pseudocyclic.

\begin{teor}\label{pseudocyc}
Let $(X, \{R_i\}_{0\leq i\leq n})$ be an association scheme, and for $x\in X$ and $1\leq i\leq n$,
let $R_i(x)=\{y\mid (x,y)\in R_i\}$. Then the following are equivalent.\\
(1). $(X, \{R_i\}_{0\leq i\leq n})$ is pseudocyclic.\\
(2). For some constant $t$, we have $v_j=t$ and $\sum_{k=1}^{n}p_{kj}^k=t-1$,
for $1\leq j\leq n$.\\
(3). $(X, {\mathcal B})$ is a $2-(v,t,t-1)$ design, where
${\mathcal B}=\{R_i(x)\mid x\in X, 1\leq i\leq n\}$.
\end{teor}

For a proof of this theorem, we refer the reader to \cite[p.~48]{bcn} and \cite[p.~84]{henkthesis}. Part (2) in
the above theorem is very useful. For example, we may use it to prove that the cyclotomic scheme in
Example~\ref{cycexamp} is pseudocyclic. The proof goes as follows. First, the nontrivial valencies of the
cyclotomic scheme are all equal to $f$. Second, by (\ref{cycparam}) and noting that $-1\in C_0$, we have
\begin{eqnarray*}
\sum_{k=1}^ep_{kj}^k&=&\sum_{k=1}^e|\{z\in C_{0}\mid 1+z\in C_{j-k}\}|\\
&=&|C_0|-1=f-1\\
\end{eqnarray*}

Pseudocyclic schemes can be used to construct strongly regular graphs and distance regular graphs of diameter 3
(\cite[p.~388]{bcn}). In view of this, it is of interest to construct pseudocyclic association schemes. The
cyclotomic schemes discussed above are examples of pseudocyclic association schemes on prime-power number of
points. Very few examples of pseudocyclic association schemes on nonprime-power number of points are currently
known (see \cite{mathon}, \cite[p.~390]{bcn} and \cite{henkthesis}). 
We will give examples of pseudocyclic association schemes on nonprime-power number of points in Section 6. More 
examples of such association schemes will be given in \cite{pseudo}.

\section{\label{pgl}The group $\PGL(2,q)$ as the subgroup of $\PGL(3,q)$ fixing a
nonsingular conic in $\PG(2,q)$}
Through the usual identification of $\bFqp$ with $\PG(1,q)$ given by
\[ x \leftrightarrow (x,1)^\top, \qquad \infty \leftrightarrow (1,0)^\top, \]
the 2-dimensional projective linear group $\PGL(2,q)$ acts on $\bFqp$, with action given by \beql{Gact} \forall
A= \left(
\begin{array}{cc}
a&b\\
c&d
\end{array}\right)\in \PGL(2,q),\;\mbox {and}\; \forall x\in \bFqp,\; A\cdot x=A(x):=\frac{ax+b}{cx+d} \eeql
For any four-tuple $(\ga, \gb, \gc,\gd)$ in $(\bFqp)^4$ with no three of $\ga,\gb,\gc,\gd$ equal, we define the
{\em cross-ratio\/} $\rho(\ga, \gb, \gc,\gd)$ by
\[\rho(\ga, \gb, \gc,\gd)=\frac{(\ga-\gc)(\gb-\gd)}{(\ga-\gd)(\gb-\gc)}, \]
with obvious interpretation if one or two of $\ga, \gb, \gc,\gd$ are equal to $\infty$. For example, if
$\ga=\infty$, then we define
\[
\rho(\infty,\gb, \gc,\gd)=
            \left\{ \begin{array}{ll}
       \frac{\gb-\gd}{\gb-\gc}, & \mbox{if $\gb,\gc,\gd\neq\infty$;} \\
       1, & \mbox{if $\gb=\infty$ (so $\gc,\gd\neq\infty$);} \\
       0, & \mbox{if $\gc=\infty$ (so $\gb,\gd\neq\infty$);} \\
       \infty, & \mbox{if $\gd=\infty$ (so $\gb,\gc\neq\infty$).}
                     \end{array}
            \right.
\]
(We will return to this interpretation later on.) Note that the cross-ratio is contained in $\bFqp$; moreover,
\beql{cr31} \rho(\ga, \gb, \gd,\gc)=\rho(\gb, \ga, \gc,\gd)=1/\rho(\ga, \gb, \gc,\gd)\eeql and \beql{cr32}
\rho(\gb, \ga, \gd,\gc)=\rho(\ga, \gb,\gc,\gd).\eeql Also, it is easily verified that \beql{cr33}
\mbox{$\rho(\ga, \gb, \gc,\gd)=1$ if and only if $\ga=\gb$ or $\gc=\gd$.} \eeql Observe that, with the above
identification of $\bFqp$ with $\PG(1,q)$, if $v_\ga=(\ga_0, \ga_1)^\top, v_\gb=(\gb_0, \gb_1)^\top,
v_\gc=(\gc_0, \gc_1)^\top$, and $v_\gd=(\gd_0, \gd_1)^\top$ are the four points in $\PG(1,q)$ corresponding to
$\ga, \gb, \gc$ and $\gd$ in $\bFqp$, respectively, then $\rho(\ga, \gb, \gc,\gd)$ can be identified with the
point \beql{cr30} ((\ga_0\gc_1-\ga_1\gc_0)(\gb_0\gd_1-\gb_1\gd_0),
(\ga_0\gd_1-\ga_1\gd_0)(\gb_0\gc_1-\gb_1\gc_0))^\top, \eeql of $\PG(1,q)$, which can be more conveniently written
as \beql{cr-proj} \left(
\begin{array}{c}
\det(v_\ga,v_\gc)\det(v_\gb,v_\gd)\\
\det(v_\ga,v_\gd)\det(v_\gb,v_\gc)
\end{array}
\right) \eeql Note that this last expression in~(\ref{cr-proj}) is equal to the zero vector only if three of the
four vectors $v_\ga, v_\gb, v_\gc, v_\gd$ are equal, which we have excluded. 
Therefore, (\ref{cr-proj}) allows us to interpret the value of the cross-ratio as an element in $\PG(1,q)$.

We will need several well-known properties concerning the above action of $\PGL(2,q)$
and its relation to the cross-ratio.
\begin{teor} \label{pgl-prop}
(i) The action of $\PGL(2,q)$ on $\bFqp$ defined in (\ref{Gact}) is sharply
3-transitive.

(ii) The group $\PGL(2,q)$ leaves the cross-ratio on $\bFqp$ invariant, that is,
if $A\in\PGL(2,q)$, then
$\rho(A(\ga),A(\gb),A(\gc),A(\gd))=\rho(\ga,\gb,\gc,\gd)$ for all $\ga,\gb,\gc, \gd\in
\bFqp$ with no three of $\ga, \gb, \gc, \gd$ equal.

(iii) Moreover, if $\Omega_+=\{\{\ga,\gb\}\mid \ga, \gb\in\bFqp, \ga\neq \gb\}$, then the action of $\PGL(2,q)$
on $\Omega_+\times \Omega_+$ has orbits
\[ O_{\rm diag} = \{(\{\ga,\gb\}, \{\ga,\gb\})\mid \{\ga,\gb\}\in\Omega_+\}, \]
and
\[ O_{\{r,r^{-1}\}}=\{(\{\ga,\gb\}, \{\gc,\gd\})\mid \{\ga,\gb\}, \{\gc,\gd\} \in
\Omega_+, \{\ga,\gb\}\neq \{\gc,\gd\}, \rho(\ga,\gb,\gc,\gd)\in\{r, r^{-1}\}\},
\]
for $r\in(\bFqp)\setminus\{1\}$.
\end{teor}
\begin{proof} (Sketch) It is easily proved that the triple $(\infty, 0,1)$ can be
mapped to any other triple $(\ga,\gb,\gc)$ with $\ga, \gb, \gc$ all distinct. So $\PGL(2,q)$ acts 3-transitively
on $\bFqp$. Since $\PGL(2,q)$ has size $(q^2-1)(q^2-q)/(q-1)=(q+1)q(q-1)$, part (i) follows.

From the representation (\ref{cr-proj}) of the cross-ratio, we immediately see that
$\PGL(2,q)$ indeed leaves the cross-ratio invariant, so part (ii) holds.

We have that $\rho(\infty, 0, 1, \delta)=\delta$ for all $\gd\in\bFqp$. Also, for $\{\ga,\gb\}, \{\gc,\gd\}
\in\Omega_+$, we have that $\rho(\ga,\gb,\gc,\gd)\in\{0,\infty\}$ if (and only if) $\{\ga,\gb\}\cap
\{\gc,\gd\}\neq \emptyset$. These observations are sufficient to conclude that $\rho$ takes on all values in
$\bFqp\setminus\{1\}$ and that the orbits are indeed as stated in part (iii).
\end{proof}

For any element $\xi$ in some extension field $\bF_{q^m}$ of $\bF_q$, we define a point $P_\xi$ in $\PG(2,q^m)$
by
\[P_\xi = (\xi, \xi^2, 1)^\top;\]
furthermore, we define
\[ P_{\infty}=(0,1,0)^\top\]
and
\[ P_{\rm Nuc}=(1,0,0)^\top.\]
We will denote by $\cO_{q^m}$ the subset of size $q^m+1$ of $\PG(2,q^m)$ consisting of the points $P_\xi$, where
$\xi \in \bF_{q^m}\cup\{\infty\}$. It is easily verified that for each $m$, the set $\cO_{q^m}$ is a nonsingular
conic in $\PG(2,q^m)$, and constitutes an oval. We will mostly write $\cO$ to denote $\cO_q$ and $\Obar$ to
denote $\cO_{q^2}$. For each $\xi\in\bF_{q^m}$, there is a unique tangent line through $P_\xi$ given by
\beql{tanxi} t_\xi= (2\xi,1, \xi^2)^\perp \eeql if $\xi\neq \infty$, and \beql{taninf}
t_\infty=(0,0,1)^\perp.\eeql Note that $t_\xi$ is contained in $\PG(2,q)$ if and only if $\xi\in \bFqp$. Also
note that if $q$ is even, then the point $P_{\rm Nuc}$ is the {\em nucleus\/} of the conic, that is, all tangent
lines to $\cO$ meet at the point $P_{\rm Nuc}$.

The group $\PGL(2,q)$ can be embedded as a subgroup $G$ of $\PGL(3,q)$ fixing $\cO$
setwise, by letting
\beql{GOact}
A=\left(
\begin{array}{cc}
a&b\\
c&d
\end{array}\right)
\mapsto
\left(
\begin{array}{ccc}
ad+bc&ac&bd\\
2ab& a^2&b^2\\
2cd&c^2&d^2
\end{array}\right).
\eeql
Indeed, we have the following.
\begin{teor} \label{PGLactO} Under the embedding (\ref{GOact}), the group $\PGL(2,q)$
fixes $\cO_{q^m}$ setwise for each $m$; in particular, an element $A\in\PGL(2,q)$ maps
a point $P_\xi$ on $O_{q^m}$ to the point $P_{A(\xi)}$, where $A(\xi)$ is defined as 
in ($\ref{Gact}$).
\end{teor}
\begin{proof}
It is easily verified that the image of $A$ (which we will again denote by
$A$) maps any point
$P_\xi=(\xi,\xi^2,1)^\top$ to
the vector $((a\xi+b)(c\xi+d), (a\xi+b)^2, (c\xi+d)^2)^\top$, which represents
the point $P_{A(\xi)}$.
So indeed $G$ fixes $\cO_{q^m}$ setwise.
\end{proof}

\noindent{\bf Remark:} If we identify $\cO_{q^m}$ with $\bF_{q^m}\cup\{\infty\}$ by letting
\[ P_\xi \leftrightarrow \xi,\]
then $G$ acts on $\cO_{q^m}$ in exactly the same way as $\PGL(2,q)$ acts on $\bF_{q^m}\cup\{\infty\}$ with the
action given in (\ref{Gact}). In fact, it turns out that $G$ is the full subgroup $G(\cO)$ of $\PGL(3,q)$ fixing
$\cO$ setwise, see e.g. \cite[p.~158]{hirsch1}. This can easily verified along the following line. Assume that a
matrix $A$ in $\PGL(3,q)$ fixes $\cO$ setwise. Then for each $x$ in $\bFqp$ the image $A P_x$ is on $\cO$, hence
satisfies the equation $X^2=YZ$. Working out this condition results in a polynomial of degree three that has all
$x\in\bF_q$ as its roots. Therefore, for $q>3$ all coefficients of the polynomial have to be zero, implying that
$A$ must have the form as described above. For $q=2,3$, the claim is easily verified directly.
\section{\label{ass}A coherent configuration containing two association schemes}
The action of the subgroup $G=G(\cO)$ of $\PGL(3,q)$ fixing the conic $\cO$ as described in the previous section
produces a coherent configuration $\cR=\cR(q)$ on the set $\cL$ of non-tangent lines of $\cO$ in $\PG(2,q)$. Here
we will determine the orbits of $G(\cO)$ on $\cL\times \cL$, and show that we obtain association schemes on both
the set $\cL_+$ of hyperbolic lines and the set $\cL_-$ of elliptic lines. First, we need some preparation.

In what follows, we will repeatedly consider ``projective objects'' over a base field as a subset of similar 
projective objects over an extension field. (For example, we will consider $\PG(2,q)$ as a subset of $\PG(2,q^2)$
and $\PGL(2,q)$ as a subset of $\PGL(2,q^2)$.) In such situations it is crucial to be able to determine whether a
given projective object over the extension field is actually an object over the base field. The next theorem
addresses this question.
\begin{teor}\label{ext}Let $\bF$ be a field and let $\bE$ be a Galois extension of
$\bF$, with Galois group $\Gal$. Let $A$ be an $n\times m$ matrix with entries from
$\bE$. Then there exists some $\lambda\in\bE\setminus\{0\}$ such that $\lambda A$ has
all its entries in $\bF$ if and only if for all $\sigma\in\Gal$ there exists some
$\mu_\sigma$ in $\bE$ such that $A^\sigma=\mu_\sigma A$.
\end{teor}
\begin{proof}
Note that given $x\in \bE$, we have $x\in\bF$ if and only if $x^\sigma=x$ for all
$\sigma\in\Gal$. \\
(i) If $\lambda\in\bE\setminus\{0\}$ such that $\lambda A$ has all its entries in
$\bF$, then for $\sigma\in\Gal$, we have $\lambda^\sigma A^\sigma=\lambda A$, hence
with $\mu_\sigma=\lambda/\lambda^\sigma$, we have that $A^\sigma=\mu_\sigma A$.\\
(ii) Conversely, suppose that $A^\sigma=\mu_\sigma A$ for every $\sigma\in\Gal$.
If $A=0$, then we can take $\lambda=1$. Otherwise, let $a$ be some nonzero entry of
$A$. Since $A^\sigma=\mu_\sigma A$, we have that $a^\sigma=\mu_\sigma a$. Set
$\lambda=a^{-1}$. Then $\lambda^\sigma=(a^{-1})^\sigma=(a^\sigma)^{-1}$, hence $\lambda
=\lambda^\sigma \mu_\sigma$. As a consequence, $(\lambda A)^\sigma=\lambda^\sigma
A^\sigma = (\lambda/\mu_\sigma) \mu_\sigma A=\lambda A$. Since this holds for all
$\sigma \in \Gal$, we conclude that $\lambda A$ has all its entries in $\bF$.
\end{proof}

\begin{remark}
The usual method to prove that some scalar multiple $\gl A$ of a matrix $A$ has all entries in the base field is
to take $\gl=a^{-1}$, for some nonzero entry $a$ of $A$. (It is easy to see that if such a scalar exists then
this choice must work.) However, this approach often requires a similar but distinct argument for each entry of
$A$ separately. The above theorem can be used to avoid such inelegant case distinction, and therefore deserves to
be better known. Although the result is unlikely to be new, we do not have a reference.
\end{remark}

Consider a point $P=(x,y,z)^\top$ in $\PG(2,q^2)$. If some nonzero multiple $\lambda P$ has all its coordinates
in $\bF_q$, then we may regard $P$ as actually belonging to $\PG(2,q)$. Let us call such points {\em real\/}, and
the remaining points in $\PG(2,q^2)$ {\em virtual\/}. Similarly, we will call a line in $\PG(2,q^2)$ {\em real\/}
if it contains at least two real points, and {\em virtual\/} otherwise. It is not difficult to see that each real
line $\ell=(a,b,c)^\perp$ in fact contains $q+1$ real points and that $\ell$ is real if and only if some nonzero
multiple $\lambda (a,b,c)$ has all its entries in $\bF_q$. As a consequence, the real points in $\PG(2,q^2)$
together with the real lines in $\PG(2,q^2)$ constitute the plane $\PG(2,q)$, a Baer subplane in $\PG(2,q^2)$.

Now let $\ell\in\cL$ be any non-tangent line to $\cO$ in $\PG(2,q)$. Then $\ell$ extends to a real line in
$\PG(2,q^2)$ (which by abuse of notation we shall again denote by $\ell$). By inspection of (\ref{tanxi}) and
(\ref{taninf}), we see that all tangent lines $t_\xi$ to $\Obar=\cO_{q^2}$ in $\PG(2,q^2)$ are either virtual
tangent lines (if $\xi \in \bF_{q^2}\setminus \bF_q$) or real tangent lines in $\PG(2,q)$ (if $\xi \in\bFqp$);
therefore $\ell$ intersects $\Obar$ in two points, $P_\alpha$ and $P_\beta$, say. In fact it is easily seen that 
either $\ell$ is hyperbolic (i.e., $\alpha, \beta \in \bFqp$), or $\ell$ is elliptic (i.e., $\beta = \alpha^q$
with $\alpha\in \bF_{q^2}\setminus\bF_q$). We will let $\cL_+$ and $\cL_-$ denote the set of hyperbolic and
elliptic lines, respectively, and we will say that a line in $\cL_+$ (respectively $\cL_-$) is of {\em hyperbolic
type\/} (respectively, of {\em elliptic type\/}). Also, we define
$$\Omega_+=\{\{\ga,\gb\}\mid \ga,\gb\in \bFqp,\; \ga\neq \gb\},\;\Omega_-=\{\{\ga,\gb\}\mid \gb=\ga^q,\;\ga\in\bF_{q^2}\setminus\bFq\},$$
and
$$\Omega=\Omega_+\cup\Omega_-.$$
Note that according to the above remarks, there is a one-to-one correspondence between lines in $\cL_\gre$ and
pairs in $\Omega_\gre$ such that $\ell\in\cL_\gre$ corresponds to $\{\ga,\gb\}\in\Omega_\gre$ if $\ell\cap
\cO_{q^2}=\{P_\ga,P_\gb\}$. Also note that if $\ell$ and $m$ are two lines in $\cL$, with corresponding pairs
$\{\ga,\gb\}$ and $\{\gc,\gd\}$ in $\Omega$, respectively, and if $g_A$ is an element of $G(\cO)$ corresponding
to $A\in \PGL(2,q)$, then $g_A$ maps $\ell$ to $m$ precisely when $A$ maps $\{\ga,\gb\}$ to $\{\gc,\gd\}$, that
is, if $\{\gc,\gd\}=\{A(\ga),A(\gb)\}$. So the action of $G(\cO)$ on $\cL$ and that of $\PGL(2,q)$ on $\Omega$
are equivalent.

\begin{defi} Let $\ell, m$ be two non-tangent lines in $\PG(2,q)$, and suppose that
$\ell\cap \Obar=\{P_\ga, P_\gb\}$ and $m\cap \Obar=\{P_\gc, P_\gd\}$. We define the {\em cross-ratio}
$\rho(\ell,m)$ of the lines $\ell$ and $m$ as $\rho(\ell,m)=\{r,r^{-1}\}$, where $r\in\bF_{q^2}\cup\{\infty\}$ is
defined by
\[ r =\rho(\ga,\gb,\gc,\gd).\]
\end{defi}

We will now show that the cross-ratio essentially determines the orbits of $G(\cO)$ on $\cL\times \cL$. The
precise result is the following:

\begin{teor}\label{orbits}
Given two ordered pairs of non-tangent lines $(\ell,m)$ and $(\ell', m')$ with $\ell\neq m$ and $\ell'\neq m'$,
there exists an element of $G(\cO)$ mapping $(\ell,m)$ to $(\ell', m')$ if and only if

\noindent (i) $\ell$ and
$\ell'$ are of the same type;

\noindent (ii) $m$ and $m'$ are of the same type; and

\noindent (iii) $\rho(\ell,m)=\rho(\ell',m')$.
\end{teor}

\begin{proof} We first show that (i), (ii) and (iii) are necessary. Let $\ga, \gb,\gc,\gd,\ga',\gb',\gc',\gd'\in\bF_{q^2}\cup\{\infty\}$ be such that

\[ \ell\cap \cO_{q^2}=\{P_\ga, P_\gb\}, \quad m\cap \cO_{q^2}=\{P_\gc, P_\gd\},
\quad \ell'\cap \cO_{q^2}=\{P_{\ga'}, P_{\gb'}\}, \quad
m'\cap \cO_{q^2}=\{P_{\gc'}, P_{\gd'}\}.\]

As already remarked above, there exists some element $g_A\in G(\cO)$ mapping $\ell$ to $\ell'$ and $m$ to $m'$ if
and only if, under the action as in (\ref{Gact}), the associated matrix $A\in\PGL(2,q)$ maps $\{\ga,\gb\}$ to
$\{\ga',\gb'\}$ and $\{\gc,\gd\}$ to $\{\gc',\gd'\}$. Now any element of $G(\cO)$ obviously maps a hyperbolic
line to a hyperbolic line and an elliptic line to an elliptic line, hence (i) and (ii) are indeed necessary; and
by Theorem~\ref{pgl-prop}, part~(ii), after interchanging $\gc'$ and $\gd'$ if necessary, we have
$\rho(\ga,\gb,\gc,\gd)=\rho(\ga',\gb',\gc',\gd')$. So we see that (iii) is also necessary.

Conversely, assume that the conditions (i), (ii) and (iii) hold. By applying Theorem~\ref{pgl-prop}, part~(iii),
with $q^2$ in place of $q$, we conclude from condition (iii) that (after interchanging $\gc'$ and $\gd'$ if
necessary) there exists a (unique) matrix $A\in\PGL(2,q^2)$  mapping $\ga$ to $\ga'$, $\gb$ to $\gb'$, $\gc$ to
$\gc'$, and $\gd$ to $\gd'$. We have to show that actually $A\in\PGL(2,q)$, that is, some nonzero multiple $\gl
A$ of $A$ has all its entries in $\bF_q$. So let
\[
A=\left(
\begin{array}{cc}
a&b\\
c&d
\end{array}\right).
\]
According to our assumptions, we first have that $A$ maps $(\ga, \gb)$ to~$(\ga',\gb')$, that is, we have
\beql{Aellellp} \frac{a\ga+b}{c\ga+d}=\ga', \qquad \frac{a\gb+b}{c\gb+d}=\gb'. \eeql We distinguish two cases.

If both $\ga, \ga'\in\bFqp$, then also $\gb, \gb'\in\bFqp$.
Now from (\ref{Aellellp}) we conclude that
\[ \frac{a^q\ga+b^q}{c^q\ga+d^q} = \frac{a\ga+b}{c\ga+d}, \]
hence $\ga$ is a zero of the polynomial
\beqa
F_A(x) &=& (a^q x +b^q)(cx+d)-(ax+b)(c^qx+d^q)\\
    &=&(a^qc-ac^q)x^2
+(a^qd-ad^q+b^qc-bc^q)x + (b^qd-bd^q). \eeqa Note that this also holds for $\ga=\infty$ if we adopt the
convention that a polynomial of degree at most two has $\infty$ as a zero if and only if the polynomial has
actually degree at most one. Indeed, $F_A$ has $\infty$ as its zero if and only if $a/c=a^q/c^q$, and
$\ga'=A(\infty)=a/c$. So we conclude that if $\ell$ and $\ell'$ are both hyperbolic, then $\ga$, and by a similar
reasoning also $\gb$, are zeroes of the polynomial $F_A(x)$.

On the other hand, if both $\ga, \ga'\in\bFqmin$, then also $\gb=\ga^q$ and
$\gb'=\ga'^q$ are in $\bFqmin$.
By raising the second equation in (\ref{Aellellp}) to the $q$-th power, we again
conclude that
\[ \frac{a^q\ga+b^q}{c^q\ga+d^q} = \frac{a\ga+b}{c\ga+d}, \]
hence again we have that $\ga$, and similarly $\ga^q$, is a zero of the polynomial
$F_A(x)$.

In summary, if $A$ maps $(\ga, \gb)$ to~$(\ga',\gb')$, we can conclude that both $\ga$ and $\gb$ are zeroes of
$F_A$; hence according to our assumptions all four of $\ga, \gb, \gc, \gd$ determined by the lines $\ell$ and $m$
are zeroes of the polynomial $F_A(x)$. Now since $\ell\neq m$, we have $|\{\ga,\gb,\gc,\gd\}|\geq 3$.
Consequently $F_A(x)$ is the zero polynomial, that is, \beql{FAeq} ac^q \in\bF_q, \qquad bd^q \in \bF_q, \qquad
a^qd-bc^q =ad^q-b^qc \in\bF_q. \eeql Now we want to apply Theorem~\ref{ext}. With
\[ \Phi=a^qd - bc^q = ad^q-b^qc, \qquad \Delta=\det(A)=ad-bc \neq 0,\]
we have that
\beqa
a\Phi &=& a(a^qd-bc^q) = a^{q+1}d-ba^qc = a^q \Delta;\\
b\Phi &=& b(ad^q-b^qc) = ab^qd-b^{q+1}c = b^q \Delta;\\
c\Phi &=& c(a^qd-bc^q) = c^qad-bc^{q+1} = c^q \Delta;\\
d\Phi &=& d(ad^q-b^qc) = ad^{q+1}-d^qbc = d^q \Delta; \eeqa hence $A^q \Delta=A \Phi$, i.e., $A^q = (\Phi/\Delta)
A$. By Theorem~\ref{ext}, we may now conclude that essentially $A\in\PGL(2,q)$.
\end{proof}

\begin{cor}\label{gentrans}
The group $G(\cO)$ is generously transitive on both $\cL_+$ and $\cL_-$.
\end{cor}
\begin{proof}
Let $\ell,m$ be two distinct lines in $\cL$. Obviously, $\rho(\ell,m)=\rho(m,\ell)$.
Hence according to Theorem~\ref{orbits}, there is an element in $G(\cO)$ that maps
$\ell$ to $m$, or, equivalently, $(\ell,m)$ to $(m,\ell)$, if and only if $\ell$ and
$m$ are of the same type.
\end{proof}

Our next result relates the value of the cross-ratio $\rho(\ell,m)$ of two lines
$\ell$ and $m$
to their types. Let us define the subsets $\bB_0$ and $\bB_1$ of
$\bF_{q^2}\cup\{\infty\}$ by
\[ \bB_0=(\bF_q\cup\{\infty\})\setminus\{1\}, \qquad
\bB_1=\{x\in\bF_{q^2}\setminus\{1\}\mid x^q=x^{-1}\}. \] Note that $|\bB_0|=|\bB_1|=q$, also the intersection of
$\bB_0$ and $\bB_1$ is empty if $q$ is even, and contains only $-1$ if $q$ is odd. We have the following.
\begin{lem}\label{type} Let $\ell, m$ be two distinct non-tangent lines in
$\PG(2,q)$, and let $\rho(\ell,m)=\{\gl,\gl^{-1}\}$, where $\gl\in \bF_{q^2}\cup\{\infty\}$. Then $\gl$ is
contained in~$\bB_0$ if $\ell$ and $m$ are of the {\em same\/} type, and contained in~$\bB_1$ if
$\ell$ and $m$ are of {\em different\/} type.
\end{lem}

\begin{proof}
Easy consequence of the fact that if $\ga, \gb, \gc, \gd\in\bFqp$ and
$\xi,\eta\in\bF_{q^2}\setminus\bF_q$, then $\rho(\ga,\gb,\gc,\gd)$ and
$\rho(\xi,\xi^q,\eta,\eta^q)$
are both in~$\bB_0$ while $\rho(\xi,\xi^q,\gc,\gd)$ and $\rho(\ga,\gb,\eta,\eta^q)$
are both in~$\bB_1$.
\end{proof}


\vskip5pt

For $\gre, \gf\in\{1,-1\}$ and for $\gl\in\bB_0$ (if $\gre=\gf=1$), or $\gl\in\bB_0\setminus\{0,\infty\}$ (if
$\gre=\gf=-1$), or  $\gl\in\bB_1$ (if $\gre\neq\gf$), we define
\[ \cR_{\{\gl,\gl^{-1}\}}(\gre,\gf)=\{(\ell,m)\in \cL_\gre\times\cL_\gf , \ell\neq m
\mid \rho(\ell,m)=\{\gl,\gl^{-1}\}\} . \] We observed earlier that $\rho(\ell,m)\neq1$ and
$\rho(\ell,m)=\{0,\infty\}$ if and only if $\ell$ and $m$ are equal or intersect on $\cO$. Hence according to
Theorem~\ref{orbits} and Lemma~\ref{type}, each of the non-diagonal orbits of $G(\cO)$ on $\cL\times\cL$, that
is, each non-diagonal relation of the coherent configuration $\cR$ obtained from the action of $G(\cO)$ on
$\cL\times\cL$, is actually of the form $\cR_{\{\gl,\gl^{-1}\}}(\gre,\gf)$ with the restrictions on $\gl$ as
given above. Moreover, since $G(\cO)$ is transitive on both $\cL_+$ and $\cL_-$, we have that
\[ |\cR_{\{\gl,\gl^{-1}\}}(\gre,\gf)| = |\cL_\gre| v_{\{\gl,\gl^{-1}\}}(\gre,\gf),\]
where the numbers $v_{\{\gl,\gl^{-1}\}}(\gre,\gf)$ are the {\em valencies\/} of the coherent configuration $\cR$.
In order to finish our description of the orbits of $G(\cO)$ on $\cL\times \cL$, we will show that each of the
orbits defined above is indeed nonempty.

\begin{teor}\label{vals}
We have that
\[
v_{\{\gl,\gl^{-1}\}}(\gre,\gf)=
            \left\{ \begin{array}{ll}
               2(q-1), & \mbox{if $\gre=\gf=1$ and $\{\gl,\gl^{-1}\}=\{0,\infty\}$;} \\
               (q-\gre)/2, & \mbox{if $q$ is odd and $\gl=-1$;} \\
               q-\gre, & \mbox{if $\gl\in\bB_{(1-\delta_{\epsilon,\phi})}$
                and $\gl\neq -1, 0, \infty$.}
             \end{array}
            \right.
\]
(Here $\delta$ is the Kronecker delta.)
\end{teor}

\begin{proof} Fix a line $\ell\in\cL_\gre$, and let $\ell\cap \cO_{q^2}=\{P_\ga,P_\gb\}$, where
$\{\ga,\gb\}\in\Omega_\gre$. We want to count the number of $m\in \cL_\gf$ such that $m \cap\cO_{q^2}=\{P_x,
P_y\}$ with $\{x,y\}\in\Omega_\gf$, and $\rho(\ell,m)=\{\gl,\gl^{-1}\}$, where \beql{rhoeq}
\lambda=\rho(\ga,\gb,x,y) = \frac{(\ga-x)(\gb-y)} {(\ga-y)(\gb-x)}.\eeql Now we note the following. First, we
have $\gl\in\{0,\infty\}$ if and only if $\{\ga,\gb\}\cap\{x,y\}\neq\emptyset$, that is, if and only if the
corresponding lines $\ell$ and $m$ intersect on $\cO$. Hence
\[v_{\{0,\infty\}}(\gre,\gf)=\left\{\begin{array}{ll} 2(q-1), &\mbox{if $\gre=\gf=1$};\\
 0, & \mbox{otherwise}.\end{array}\right.\] Next, we have $x=y$ in
(\ref{rhoeq}) only if $\gl=1$, which is excluded. Also, by interchanging $x$ and $y$ the cross-ratio $\gl$ in
(\ref{rhoeq}) is {\em inverted\/}, and the only cases where $\gl=\gl^{-1}$ are $\gl=1$ (which is excluded) and
$\gl=-1$. As a consequence, for $\gl\in(\bB_0\cup\bB_1)\setminus\{0,\infty\}$ the number
$v_{\{\gl,\gl^{-1}\}}(\gre,\gf)$ equals the number of solutions $(x,y)$ of (\ref{rhoeq}) with
$\{x,y\}\in\Omega_\gf$ if $q$ is even or $\gl\neq -1$, and is equal to half of the number of such solutions if
$\gl=-1$ and $q$ is odd.

First, let $\gre=1$. According to Theorem~\ref{orbits}, we may assume without loss of generality that
$\ga=\infty$ and $\gb=0$, so that (\ref{rhoeq}) reduces to $\gl=-y/(-x)=y/x$. If $\gf=1$, then
$x,y\in\bF_q\cup\{\infty\}$; so evidently we must have $\gl\in\bF_q\setminus\{0,1\}$; in that case for each
$x\in\bF_q\setminus\{0\}$ there is a unique solution $y\in\bF_q$, so there are $q-1$ solutions in total. 
Similarly, if $\gf=-1$, then $x\in\bFqmin$ and $y=x^q$, so (\ref{rhoeq}) reduces to $\gl=x^{q-1}$. 
Hence $\gl\in\bB_1$ and again there are precisely $q-1$ solutions for each such $\gl$.

If $\gre=-1$, then we have $\ga\in\bFqmin$ and $\gb=\ga^q$. First, if $\gf=1$, then $x,y\in\bFqp$. In that case
we see immediately from (\ref{rhoeq}) that $\lambda^q=1/\lambda$, hence there are no solutions except when
$\lambda\in \bB_1$. Let $z$ be the unique solution of the equation $\lambda = (\ga-z)/(\ga^q-z)$. Now for
$y=\infty$ the unique solution is $x=z$; for $y=z$ the unique solution is $x=\infty$, and it is easily seen that
for each $y\in\bF_q\setminus\{z\}$ there is a unique solution $x\in\bF_q$. So there are precisely $q+1$
solutions. Finally, if $\gf=-1$, then we have $x\in\bFqmin$ and $y=x^q$. In that case, the desired solutions of
(\ref{rhoeq}) satisfy
\[ (\ga-x)(\ga^q-x^q) = \gl (\ga-x^q)(\ga^q-x),\] with $x\in\bFqmin$ and $x\neq\ga, \ga^q$.
Clearly this is only possible if $\gl\in\bFq\setminus\{0,1\}$. For each such $\gl$, this equation has at most
$q+1$ solutions in $\bFqmin$. On the other hand, there are $q^2-q-2$ choices of $x$ with $x\neq \ga,\ga^q$;
consequently, the {\em average\/} number of valid solutions equals $(q^2-q-2)/(q-2)=q+1$. Since the average
number of solutions equals the maximum number of solutions, there must be exactly $q+1$ solutions for each
$\gl\in\bF_q\setminus\{0,1\}$, and the result stated for $\gre=\gf=-1$ follows.
\end{proof}

By combining Theorem~\ref{orbits} and Theorem~\ref{vals} we obtain the following result.

\begin{teor}\label{main-res}
(i) The action of the group $G(\cO)$ on $\cL$ affords a weakly symmetric coherent configuration $\cR(q)$.

\noindent (ii) The restriction of $\cR(q)$ to the the fibre $\cL_+$ of hyperbolic lines constitutes an
association scheme $\cR_{+}(q)$ (or $\cH(q)$), with $q/2$ classes if $q$ is even and with $(q+1)/2$ classes if
$q$ is odd. The non-diagonal relations are precisely the sets $\cR_{\{\gl,\gl^{-1}\}}(1,1)$ with $\gl\in \bB_0$,
with corresponding valencies $v_{\{\gl,\gl^{-1}\}}=v_{\{\gl,\gl^{-1}\}}(1,1)$.

\noindent (iii) The restriction of $\cR(q)$ to the the fibre $\cL_-$ of elliptic lines constitutes an
association scheme $\cR_{-}(q)$ (or $\cE(q)$), with $q/2-1$ classes if $q$ is even and with $(q-1)/2$ classes if
$q$ is odd. The non-diagonal relations are precisely the sets $\cR_{\{\gl,\gl^{-1}\}}(-1,-1)$ with
$\gl\in\bF_q\setminus\{0,1\}$, with corresponding valencies $v_{\{\gl,\gl^{-1}\}}=v_{\{\gl,\gl^{-1}\}}(-1,-1)$.
\end{teor}

We will refer to the association schemes $\cH(q)$ and $\cE(q)$ in part (ii) and (iii) of the above
theorem as the {\em hyperbolic\/} and {\em elliptic\/} scheme, respectively. The hyperbolic scheme was recently 
investigated in \cite{CD} as a refinement (fission) of the triangular scheme. The elliptic scheme was first described 
in \cite{henkthesis} but our approach here is new.

\section{\label{rhohat} An expression to determine the orbit from the homogeneous coordinates of the lines}
Let $\epsilon,\phi\in\{-1,1\}$ and let $(\ell,m)\in \cL_{\epsilon}\times\cL_{\phi}$. In this section we will
develop an expression $\rhoh(\ell,m)$ that can be used to index the relation of $\cR(q)$ containing $(\ell,m)$,
in terms of the homogeneous coordinates of $\ell$ and $m$.

We need some preparation. Consider the function $f: \bF_{q^2}\cup\{\infty\}\rightarrow \bF_{q^2}\cup\{\infty\}$ defined by
\[ f(x) =
            \left\{ \begin{array}{ll}
                    \frac{1}{x+x^{-1}}, & \mbox{if $q$ is even;} \\
                    \frac{1}{4}+\frac{1}{-2+x+x^{-1}}, & \mbox{if $q$ is odd,}
                     \end{array}
            \right.
\]
for $x\in\bF_{q^2}\setminus\{0,1\}$, $f(1)=\infty$, and $f(0)=f(\infty)=0$ if $q$ is even and
$f(0)=f(\infty)=1/4$ if $q$ is odd. (Note that the values of $f$ on $\infty, 0,1$ are consistent with the general
expression for $f(x)$ when we interpret $1/0=\infty$ and handle $\infty$ in the usual way.) This function has a
few remarkable properties. To describe these, we introduce some notation. For $q=2^r$ and for $e\in\bF_2$, let
$\bT_e=\bT_e(q)$ denote the collection of elements with absolute trace $e$ in $\bF_q$, that is,
\[\bT_e=\{x\in\bF_q\mid \Tr(x):=x+x^2 + \cdots +x^{2^{r-1}}=e\}\]
For $q$ odd, we let $\bT_0$ and $\bT_1$ denote the
collection of nonzero squares and non-squares in $\bF_q$, respectively, that is,
\[\bT_0=\{ x^2 \mid x \in\bF_q\} \setminus \{0\}\]
and $\bT_1 = \bF_q \setminus (\{0\}\cup \bT_0)$. Note that in the case where $q$ is even, it is well-known that
\[\bT_0=\{x^2+x \mid x\in\bF_q\}\]
\begin{lem}\label{f-prop} The function $f$ has the following properties:\\
(i) $f(x)=f(y)$ if and only if $x=y$ or $x=y^{-1}$;\\
(ii) $f(x)=\infty$ if and only if $x=1$;\\
(iii) if $q$ is odd, then $f(x)=0$ if and only if $x=-1$;\\
(iv) $f(x) \in\bF_q$ if and only if $x\in\bB_0\cup \bB_1$;\\
(v) if $q$ is even and $x\in \bF_{q^2}\setminus\{0,1\}$, then
\[ f(x) = \frac{1}{x+1}+\frac{1}{(x+1)^2},\]
and if $q$ is odd and $x\in \bF_{q^2}\setminus\{0,1,-1\}$, then
\[ f(x) = \left(\frac{x+2+x^{-1}}{2(x-x^{-1})}\right)^2.
\]
Hence for $x\in \bF_{q^2}\cup\{\infty\}$ and $e\in\bF_2$, we have that
$f(x) \in\bT_e$ if and only if $x\in\bB_e\setminus\{-1\}$.
\end{lem}
\begin{proof}
Note first that $f(x)=f(y)$ if and only if $x+1/x=y+1/y$; hence part (i) follows.
Parts (ii) and (iii) are evident. To see (iv), first note that $f(x)^q=f(x^q)$, then use
part (i) to conclude that 
$f(x)\in\bFqp$ if and only if $x^q\in\{x,x^{-1}\}$. The expressions for $f(x)$ in part (v) are easily verified.
Since $f(\infty)=0\in\bT_0$ if $q$ is even and $f(\infty)=1/2^2\in\bT_0$ if $q$ is odd, the expressions in (v)
imply that $f(x)\in\bT_0$ if and only if $x\in(\bFqp)\setminus\{1,-1\}$. Now the remainder of part (v) follows
from (iv).
\end{proof}

Next we determine the type of a line in terms of its homogeneous coordinates,
and we establish relations between the homogeneous coordinates of a line and the points
of intersection of this line with the conic $\cO_{q^2}$.

\begin{lem}\label{hom-line-rel} Let $\ell$ be a line in $\PG(2,q)$ with homogeneous coordinates
$\ell=(z,x,y)^\perp$, and let $\ell\cap \cO_{q^2}=\{P_\ga,P_\gb\}$,
where $\ga, \gb\in \bF_{q^2}\cup\{\infty\}$ and $\ga=\gb$ if $\ell$ is a tangent line. Define $\Delta(\ell)\in\bFqp$ by
\[
\Delta(\ell) =
            \left\{ \begin{array}{ll}
                    xy/z^2, & \mbox{if $q$ is even;} \\
                    1/(z^2-4xy), & \mbox{if $q$ is odd.}
                     \end{array}
            \right.
\]
(i) We have that $\ell\in\cL_{(-1)^e}$ if and only if $\Delta(\ell)\in\bT_e$, and $\ell$ is a tangent line to
$\cO$ in $\PG(2,q)$ if and only if $\Delta(\ell)=\infty$.

\noindent (ii) If $x\neq 0$, then \beql{gab-rel} \ga+\gb = -z/x, \qquad \ga \gb = y/x; \eeql and if $x=0$, then
\beql{gab-rel-0} \ga = \infty, \qquad \gb = -y/z. \eeql
\end{lem}
\begin{proof}
Note first that by definition $\ga, \gb$ are the solutions in $\bF_{q^2}\cup\{\infty\}$ of the quadratic equation
\[ z\xi + x\xi^2 + y =0. \]
(Here, by convention, $\xi=\infty$ is a solution if and only if $x=0$.)
Now (i) follows from the standard theory on solutions of quadratic equations and (ii)
follows from this equation by writing it in the form $x(\xi-\ga)(\xi-\gb)=0$.
\end{proof}

Now let $(\ell,m)\in\cL_\gre\times \cL_\gf$ be a pair of distinct non-tangent lines in
$\PG(2,q)$, and let $\ga,\gb,\gc,\gd$ be such that
\[\ell \cap \cO_{q^2} = \{P_\ga, P_\gb\}, \qquad m \cap \cO_{q^2} = \{P_\gc, P_\gd\}.\]
Furthermore, let $\ell$ and $m$ have homogeneous coordinates
\[ \ell=(z,x,y)^\perp, \qquad m = (\zb, \xb, \yb)^\perp. \]
In the previous section we have seen that the orbit of the action of $G(\cO)$ on $\cL\times \cL$ containing the
pair $(\ell,m)$ is $\cR_{\{\rho,\rho^{-1}\}}(\gre,\gf)$, where
\[\rho=\rho(\ga,\gb,\gc,\gd).\]
Now we define the {\em modified cross-ratio\/} $\rhoh(\ell,m)$ of the lines $\ell$ and
$m$ by
\[ \rhoh(\ell,m)=f(\rho)=
            \left\{ \begin{array}{ll}
                    \frac{1}{\rho+\rho^{-1}}, & \mbox{if $q$ is even;} \\
                    \frac{1}{4}+\frac{1}{-2+\rho+\rho^{-1}}, & \mbox{if $q$ is odd.}
                     \end{array}
            \right.
\]
We will now use the previous lemma to express $\rhoh(\ell,m)$ in terms of the homogeneous
coordinates of $\ell$ and $m$.
Let $\sigma: \bF_q\rightarrow \bT_0$ be defined by
\[ \sigma(x)=
            \left\{ \begin{array}{ll}
                     x^2+x & \mbox{if $q$ is even;} \\
                     x^2 & \mbox{if $q$ is odd.}
                     \end{array}
            \right.
\]
Then the result is as follows.
\begin{teor} \label{rhoh-expr}
If $\ell=(z,x,y)^\perp$ and $m=(\zb,\xb,\yb)^\perp$ are two non-tangent lines and if
$\Delta=\Delta(\ell)$ and $\Db=\Delta(m)$,
then
\[ \rhoh(\ell,m) =
            \left\{ \begin{array}{ll}
                  \frac{(x\yb + \xb y)^2+(x \zb + \xb z)(y \zb + \yb z)}{z^2\zb^2}
           =\sigma((x \yb+\xb y)/(z \zb)) +
             \Delta+\Db, & \mbox{if $q$ is even;} \\
                    (2 x\yb + 2 \xb y -z\zb)^2\Delta\Db/4=
            \sigma(x\yb + \xb y -\frac {z\zb} {2})\Delta\Db,
                            & \mbox{if $q$ is odd.}
                     \end{array}
            \right.
\]
\end{teor}
\begin{proof}
Let $\rho=\rho(\ga,\gb,\gc,\gd)$ with $\ga, \gb, \gc$ and $\gd$ as given above, and let $\rhoh=f(\rho)$.
Initially, we will assume that $x, \xb\neq 0$. First we observe that \beql{rho-tn-exp} -2 + \rho + \rho^{-1} =
    \frac{(\ga-\gb)^2(\gc-\gd)^2}{(\ga-\gc)(\ga-\gd)(\gb-\gc)(\gb-\gd)}.
\eeql Now $(\ga-\gb)^2 = (\ga+\gb)^2 -4\ga \gb = (-z/x)^2-4y/x$, and similarly $(\gc-\gd)^2 =
(-\zb/\xb)^2-4\yb/\xb$, hence \beql{rho-n} (\ga-\gb)^2 (\gc-\gd)^2 =
            \left\{ \begin{array}{ll}
                  z^2\zb^2/(x^2\xb^2), &  \mbox{if $q$ is even;} \\
                    1/(\Delta\Db x^2\xb^2), & \mbox{if $q$ is odd.}
                     \end{array}
            \right.
\eeql
Moreover, straightforward but somewhat tedious computations show that \beqa
  &&(\ga-\gc)(\ga-\gd)(\gb-\gc)(\gb-\gd)\\
    &=& (\ga\gb)^2 - \ga\gb(\ga+\gb)(\gc+\gd) + (\ga+\gb)^2 \gc\gd +
        \ga\gb(\gc^2+\gd^2)-(\ga+\gb)\gc\gd(\gc+\gd)+(\gc\gd)^2\\
    &=& ((x \yb - \xb y)^2 + (x \zb - \xb z)(y \zb - \yb z))/(x^2\xb^2).
\eeqa

By combining these expressions we obtain in a straightforward way the desired expressions for $\rhoh$. Finally,
it is not difficult to check that the expressions for $\rhoh$ are also correct in the case where one of $x, \xb$
is equal to zero.
\end{proof}

In what follows, we will use the elements of $\bF_q$ to index the relations of the coherent configuration
$\cR=\cR(q)$, and the modified cross-ratio $\rhoh$ to determine the relation of a given pair of distinct
non-tangent lines. For $\gre,\gf\in \{-1,1\}$ and $\lambda\in\bF_q$, we define
\[ \cR_{\mbox{diag}}(\gre,\gre) := \{(\ell,\ell) \mid \ell\in\cL_\gre\}
\]
and
\[ \cR_{\lambda}(\gre,\gf):=\{(\ell,m)\in\cL_\gre\times \cL_\gf \mid \ell\neq m,\;
\rhoh(\ell,m)={\gl}\}.
\]
Here $\cR_{\mbox{diag}}(1,1)$ and $\cR_{\mbox{diag}}(-1,-1)$ are the two diagonal relations on the fibres $\cL_+$
and $\cL_{-}$. Lemma~\ref{hom-line-rel} together with the expressions for $\rhoh(\ell,m)$ in
Theorem~\ref{rhoh-expr} show that the types of $\ell$ and $m$ alone determine whether $\rhoh(\ell,m)$ is
contained in $\bT_0$ or in $\bT_1$ (except in the case where $\rhoh(\ell,m)=0$ if $q$ is odd). In order to state
the next theorem concisely, we define for $e\in\bF_2$,
\[
\bT_e^+=
            \left\{ \begin{array}{ll}
                    \bT_e, & \mbox{if $q$ is even,} \\
                    \bT_e\cup \{0\}, & \mbox{if $q$ is odd;}
                     \end{array}
            \right.
\]
and $\bT_0^* = \bT_0\setminus\{0\}$.

Now a careful inspection of Theorem~\ref{main-res} in fact shows that we have the following.

\begin{teor}\label{Rnonempty}
The non-diagonal relations $\cR_{{\gl}}(\gre,\gf)$ of the coherent configuration $\cR(q)$ are nonempty precisely
when

\noindent
(i) $\gre=\gf=1$ and $\gl\in \bT_0^+$;\\
(ii) $\gre\neq\gf$ and $\gl\in\bT_1^+$; or\\
(iii) $\gre=\gf=-1$ and $\gl\in \bT_0^*$ if $q$ is even or $\gl\in \bT_0^+\setminus\{1/4\}$ if $q$ is odd.
\end{teor}

\section{The intersection parameters of $\cR(q)$ in the case of even characteristic}
{\em In the rest of this paper, we always assume that $q$ is even.\/} Here we will determine the intersection
parameters of the coherent configuration $\cR(q)$ in the case of even characteristic. In this case, the results
from the previous section can be resumed as follows. By Lemma~\ref{hom-line-rel}, a line $\ell$ in $\PG(2,q)$ is
non-tangent to $\cO$ if and only if it can be represented in homogeneous coordinates as $\ell=(1,x,y)^\perp$ with
$x,y\in\bF_q$; if $\Delta=xy\in\bT_e$, then $\ell\in\cL_\gre$ with $\gre=(-1)^e$. Moreover,
Theorem~\ref{rhoh-expr} implies that if $\ell=(1,x,y)^\perp\in\cL_\gre$ and $m=(1,z,u)^\perp\in\cL_\gf$ are two
non-tangent lines with $xy\in\bT_e$ and $zu\in\bT_f$, then $\gre=(-1)^e$, $\gf=(-1)^f$, and \beql{rhohlm}
\rhoh(\ell,m)=(xu +yz)^2 + (xu+yz) + xy + zu = x^2u^2+y^2z^2 +(x+z)(y+u). \eeql Since $(xu +yz)^2 +
(xu+yz)\in\bT_0$, we see that $\rhoh(\ell,m)$ is contained in $\bT_{e+f}$. Furthermore, we recall that
$\rhoh(\ell,m)=0$ precisely when $\ell=m$ or when $\ell$ and $m$ are lines in $\cL_1$ that intersect on $\cO$.
For later reference, we state these observations explicitly.

\begin{lem}\label{rhoh-type} Let $c\in\bF_q$, and $e,f\in\bF_2$. If $(\ell,m)\in \cR_c(\gre, \gf)$ with $\gre=(-1)^e$ and
$\gf=(-1)^f$, then $c\in \bT_{e+f}$. Moreover, if $c=0$, then $\gre=\gf=1$.
\end{lem}

\begin{cor}\label{int-type}
Let $a,b,c\in\bF_q$, and let $e,f,g\in\bF_2$. Write $\gre=(-1)^e$,
$\gf=(-1)^f$, and $\theta=(-1)^g$. If $\ell\in\cL_\gre$, $m\in\cL_\gf$, and $n\in
\cL_\theta$ with $\rhoh(\ell,m)=c$, $\rhoh(\ell,n)=a$, and $\rhoh(n,m)=b$, then
$\Tr(c)=e+f$, $\Tr(a)=e+g$, and $\Tr(b)=f+g$; in particular, $a+b+c\in\bT_0$.
\end{cor}

In order to determine the intersection parameters of $\cR(q)$, we need to do the
following. Choose any pair $(\ell,m)\in\cR_c(\gre,\gf)$, with $\gre=(-1)^e$ and
$\gf=(-1)^f$, say, and then count the number of lines $n\in \cL_\theta$, with
$\theta=(-1)^g$, say, such that $(\ell,n)\in\cR_a(\gre,\theta)$ and
$(n,m)\in\cR_b(\theta,\gf)$. Now note that according to Lemma~\ref{rhoh-type}, there
are such lines $n$ only if $\Tr(c)=e+f$, $\Tr(a)=e+g$, and $\Tr(b)=g+f$, and so, in
particular, only if $a+b+c\in\bT_0$. These observations motivate the following
definitions.

For all $a,b,c\in\bF_q$, for all $\gre,\gf\in\{-1,1\}$, and for all lines $\ell,m$ with
$(\ell,m)\in\cR_c(\gre,\gf)$ (so $\ell\neq m$), we define \beql{va-def} v_a(\ell)=v_a(\gre)=\#
\{n\in\cL\setminus\{\ell\} \mid \rhoh(\ell,n)=a\}, \eeql \beql{pabc-def} p_{a,b}(\ell,m)=p_{a,b}^c(\gre)= \#
\{n\in\cL\setminus\{\ell,m\} \mid \mbox{$\rhoh(\ell,n)=a$ and $\rhoh(n,m)=b$}\}, \eeql and \beql{piabc-def}
\pi_{a,b}(\ell,m)=\pi_{a,b}^c(\gre)= \# \{n\in\cL \mid \mbox{$\rhoh(\ell,n)=a$ and $\rhoh(n,m)=b$}\}. \eeql
The above observations show that the numbers in (\ref{va-def}) and (\ref{pabc-def}) are valencies and
intersection parameters of $\cR(q)$.

We also note that since $\rhoh(\ell,\ell)=0$ for all lines $\ell\in\cL$, we have the following:

\begin{lem}\label{pi-p}
Let $a,b,c\in\bF_q$ and let $\gre\in\{-1,1\}$. Then
\[
\pi^c_{a,b}(\gre) = p^c_{a,b}(\gre)
    + \delta_{a,0}\delta_{b,c} + \delta_{b,0}\delta_{a,c}.\]
\end{lem}

In what follows, we will sometimes use the symbol $\infty$ to indicate a diagonal relation and write
$\cR_\infty(\gre,\gre)$ instead of $\cR_{\rm diag}(\gre,\gre)$. Remark that the intersection parameters involving
a diagonal relation are $p_{a,b}^\infty(\gre)=\delta_{a,b}v_a(\gre)$, $p^c_{\infty,b}(\gre)= \delta_{b,c}$, and
$p^c_{a,\infty}(\gre)=\delta_{a,c}$.

According to Lemma~\ref{pi-p}, in order to obtain all intersection parameters, it is sufficient to compute the
numbers $\pi^c_{a,b}(\gre)$ for $a,b,c\in\bF_q$ with $a+b+c\in\bT_0$. We begin with the following observation.

\begin{lem}\label{pi-sym}
For all $a,b,c\in\bF_q$, we have $\pi^c_{a,b}(\gre) = \pi^c_{b,a}((-1)^{\Tr(c)}\gre)$
and $p^c_{a,b}(\gre) = p^c_{b,a}((-1)^{\Tr(c)}\gre)$.
\end{lem}
\begin{proof}
The number $\pi^c_{a,b}(\gre)$, with $\gre=(-1)^e$,
counts the number of lines $n\in\cL$ such that $\rhoh(\ell,n)=a$
and $\rhoh(n,m)=b$, for some pair of distinct non-tangent lines $\ell,m$ with
$\ell\in\cL_\gre$ and $\rhoh(\ell,m)=c$.
By Corollary~\ref{int-type}, we then have
$m\in\cL_\gf$ with $\gf=(-1)^f$ and  $f=e+\Tr(c)$;
hence all these lines $n$, and no others, contribute to
the number $\pi^c_{b,a}((-1)^{e+\Tr(c)})$.
This proves the first equality; the other equality follows from
Lemma~\ref{pi-p}.
\end{proof}

For $v\in\bF_q$, we define $\ell_v=(1,v,v)^\perp$.
Note that
\[\rhoh(\ell_v,\ell_{v+c})=c^2, \qquad \Delta(\ell_v)=v^2, \]
for any $c\in\bF_q$.

\begin{teor}\label{pi-nz}
Let $a,b,c\in\bF_q$ with $a+b+c\in\bT_0$ and $c\neq 0$. For each
$e\in\bF_2$ and for each $v\in\bF_q$ with $\Tr(v)=e$, we have
\beqa  \pi^c_{a,b}(\epsilon)
    &=&\sum_\tau |\{z\in\bF_q\cup\{\infty\} \mid z^2+z =v+ac/\tau^2\}| \\
        &=&\left\{ \begin{array}{ll}
               1+2 |\bT_e \cap\{ac\}|, & \mbox{if $a+b+c=0$};\\
        2\sum_\tau
            | \bT_e \cap\{ac/\tau^2\}|, &\mbox{if $a+b+c\in \bT_0^*$},
            \end{array}
    \right.
\eeqa
where $\epsilon=(-1)^e$, and the summations are over the two elements $\tau\in\bF_q$ such that $\tau^2+\tau = a+b+c$.
(If $\tau=0$, then $z=\infty$ is supposed to be the only solution of the equation $z^2+z=v+ac/\tau^2$.)
\end{teor}

\begin{proof}
We choose $\ell=\ell_{v^{1/2}}$ and
$m=\ell_{(v+c)^{1/2}}$ in the definition (\ref{piabc-def}) of  $\pi^c_{a,b}(\gre)$.
According to (\ref{rhohlm}), we conclude that $\pi^c_{a,b}(\gre)$ equals the number of
lines $n=(1,x,y)^\perp$ for which the following holds.
\beql{xyveq}
\left\{ \begin{array}{ll}
                     \rhoh(\ell_{v^{1/2}}, n) &=(v^{1/2} (x+y))^2 + v^{1/2} (x+y) + xy + v = a \\
                     \rhoh(n,\ell_{(v+c)^{1/2}}) &= ((v+c)^{1/2} (x+y))^2 + (v+c)^{1/2} (x+y) + xy + v + c = b
                     \end{array}
            \right.
\eeql
By adding the two equations, we see that $(x,y)\in\bF_q^2$ must satisfy
\[ (c^{1/2} (x+y))^2 + c^{1/2} (x+y) = a+b+c. \]
Since $a+b+c\in\bT_0$, there are two elements $\tau\in \bF_q$ such that $\tau^2+\tau=a+b+c$. Noting that $c\neq
0$ by assumption, we see that $(x,y)\in\bF_q^2$ is a solution of (\ref{xyveq}) if and only if \beql{xyeq} \left\{
\begin{array}{ll}
                     x+y &= \tau/c^{1/2}\\
                     x^2 + x\tau/c^{1/2} + \tau v^{1/2}/c^{1/2}+v\tau^2/c &=a+v
                     \end{array}
            \right.
\eeql Now we distinguish two cases. If $\tau=0$ (which is possible if and only if $a+b+c=0$), then (\ref{xyeq})
reduces to $x^2 = a+v$, which has a unique solution. Otherwise, $\tau\neq 0$, then the substitution
$z=xc^{1/2}/\tau +v^{1/2}c^{1/2}/\tau^{1/2}$ transforms (\ref{xyeq}) into the equation
\[ z^2+z=v+ac/\tau^2\]
with $z\in\bF_q$. The two cases can be conveniently combined by interpreting $z=\infty$ as the only solution of
the above equation when $\tau=0$.

To obtain the last expression, note that if $z$ runs through $\bF_q$, then $z^2+z+v$
runs through $\bT_e$ twice.
\end{proof}

\begin{cor}\label{pi-nzz}Let $b,c\in\bF_q$ with $b+c\in\bT_0$ and $c\neq 0$.
Then for each $e\in\bF_2$ we have that
\beq  \pi^c_{0,b}(\gre) =
        \left\{ \begin{array}{ll}
               1+2 \delta_{\gre,1}, & \mbox{if $b=c$};\\
        4 \delta_{\gre,1}, &\mbox{if $b+c\in \bT_0^*$};
            \end{array}
    \right.
\eeq
where $\epsilon=(-1)^e$.
\end{cor}

To complete the determination of the intersection parameters, we compute the numbers $\pi^0_{a,b}(\gre)$. (These
numbers can also be derived from knowledge of the valencies and the other intersection numbers by using relations
between these numbers that are valid in any coherent configuration, but the direct approach is much simpler.)
\begin{teor}\label{pi-z} Let $a,b\in\bF_q$. Then $\pi_{a,b}^0(-1)=0$ and
\[ \pi_{a,b}^0(1)=
    \left\{ \begin{array}{ll}
               q+1, & \mbox{if $a=b=0$};\\
               1, & \mbox{if $a=b\neq0$};\\
               2, &\mbox{if $a+b\in\bT_0^*$}.
            \end{array}
    \right.
\]
\end{teor}
\begin{proof}
It follows from Corollary~\ref{int-type}
that the numbers
$\pi_{a,b}^0(\gre)$ are nonzero only if $\gre=1$ and $\Tr(a)=\Tr(b)$.
Take $\ell= (1,0,0)^\perp$ and $m=(1,0,1)^\perp$. Note that $P_\infty=(0,1,0)$ and
$P_0=(0,0,1)$ are on $\ell$ and $P_\infty$ and $P_1=(1,1,1)$ are on $m$, hence $\ell$
and $m$ intersect on $\cO$ so indeed $\rhoh(\ell,m)=0$.
Now count the number of non-tangent lines $n=(1,x,y)^\perp$ such that
\[ \mbox{$\rhoh(\ell,n)=xy = a$ and $\rhoh(m,n)= x^2+x+xy =b$}, \]
or, equivalently,
\[ xy = a, \qquad x^2+x=a+b.\]
Now if $a+b\in\bT_0$, then the equation $x^2+x=a+b$ has two solutions $x$ in $\bF_q$, and for
each $x\neq 0$ the first equation $xy=a$ has the unique solution $y=a/x$. Finally, for
$x=0$ (which can occur only if $a=b$), there is no solution $y$ if $a\neq0$ and there
are $q$ solutions $y$ if $a=0$.
\end{proof}

By combining Lemma~\ref{pi-p} with our results for the numbers $\pi^c_{a,b}(\gre)$
we obtain all intersection parameters. For the sake of completeness, we also state the values of the valencies.
From Theorem~\ref{vals}, we obtain the following.

\begin{teor}\label{pi-vals} For $a\in\bF_q$, we have
\[ v_a(\gre) =
    \left\{ \begin{array}{ll}
               2(q-1)\delta_{\gre,1}, &\mbox{if $a=0$}; \\
               q-\gre, &\mbox{if $a\in\bF_q^*$}.
            \end{array}
    \right.
\]
\end{teor}

In even characteristic the elliptic association scheme $\cE(q)$ has all valencies equal to $q+1$. We will use the
expressions that we have derived for the intersection parameters to show that, in fact, the elliptic scheme is
pseudocyclic.
\begin{teor}\label{ell-pseudocyc}
(i) For all $b\in \bT_0^*$ and $\epsilon\in\{-1,1\}$, we have that
\[ \sum_{a\in \bT_0^*} p^a_{a,b}(\gre) = q -4\delta_{\gre,1}. \]
(ii) The elliptic association scheme $\cE(q)$ is pseudocyclic.
\end{teor}
\begin{proof}
(i) Let $\gre=(-1)^e$. Using Lemma~\ref{pi-p} and Theorems \ref{pi-z} and~\ref{pi-nz},
we see that for $b\in \bT_0^*$, we have
\beqa
\sum_{a\in \bT_0^*} p^a_{a,b}(\gre) &=& \sum_{a\in \bT_0^*} \pi^a_{a,b}(\gre)\\
    &=& -4\delta_{e,0}
    + \sum_{a\in \bT_0} 2 \sum_{\{\tau| \tau^2+\tau=b\}} |\tau^2 \bT_e\cap \{a^2\}| \\
      &=& -4\delta_{e,0} + 2 \sum_{\{\tau| \tau^2+\tau=b\}} |\tau^2 \bT_e \cap \bT_0|\\
    &=& -4\delta_{e,0} + 2 . 2 . q/4 \\
    &=& q -4\delta_{e,0}.
\eeqa
(ii) The non-diagonal relations of the elliptic scheme $\cE(q)$ are $\cR_c(-1,-1)$
with $c\in \bT_0^*$.
So the claim follows from part (i) together with Theorem~\ref{pseudocyc}.
\end{proof}

Let $q=2^r$ and let $k$ be an integer with $\gcd(k,r)=1$. The field automorphism $\tau_k: x\mapsto x^{2^k}$
provides in a natural way a fusion of the coherent configuration $\cR(q)$ as follows. Let the orbits of $\tau_k$
on $\bF_q$ be $C_0=\{0\}, C_1, \ldots, C_n$. Define new relations
\[ R^{\rm fus}_j(\gre,\gf) = \cup_{c\in C_j} \cR_c(\gre,\gf). \]
Now since
\[ p^c_{a,b}(\gre) = p^{\tau_k(c)}_{\tau_k(a),\tau_k(b)}(\gre),\]
it follows immediately that the fusion
\[ \cR^{\rm fus} = \{ R^{\rm fus}_j(\gre,\gf) | j=0,\ldots, n, \gre,\gf\in\{-1,1\}\} \]
is a coherent configuration, which obviously is again weakly symmetric. The valencies of this new coherent
configuration are of the form
\[ v^{\rm fus}_j(\gre) = \sum_{c\in C_j} v_c(\gre).\]
Now consider two association schemes obtained from this configuration by restricting to the set of hyperbolic or
elliptic lines, respectively. It is not difficult to see that such a scheme has all valencies equal precisely in 
the elliptic case with $q=2^r$ and $r$ prime. An old conjecture from \cite{henkthesis} states that in this case 
the scheme is again pseudocyclic. In a subsequent paper \cite{pseudo} we will prove this conjecture.

\section{Fusion schemes in the case of even characteristic}

In this section, we will assume that the field size is of the form $q^2$ with $q$ even. Our aim is to show that a
certain fusion of the coherent configuration $\cR(q^2)$ on the set $\cL(q^2)$ of non-tangent lines in
$\PG(2,q^2)$ is again a coherent configuration. We will write $\bS_0$ to denote the collection of elements with
absolute trace zero in $\bF_q$, that is,
\[ \bS_0=\bT_0(q); \]
and let
\[\bS_1=\bF_q\setminus \bS_0 \]
denote the collection of elements with absolute trace one in $\bFq$.
Also, we will write $\bT_0=\bT_0(q^2)$ and $\bT_1$ to denote the elements in $\bFqt$
with absolute trace equal to 0 or 1, respectively.
As before, for any set $U$, we write $U^*$ to denote the set $U\setminus\{0\}$ and
$\tilde{U}$ to denote the set $U\setminus\{0,1\}$.

We now define the following relations for two {\em distinct\/} lines $\ell$ and $m$ in
$\cL(q^2)$.
\begin{itemize}
\item $R_1: \qquad \rhoh(\ell,m)\in \bS_0^*$;
\item $R_2: \qquad \rhoh(\ell,m)\in \bS_1$;
\item $R_3: \qquad \rhoh(\ell,m)\in \bT_0\setminus \bF_q$;
\item $R_4: \qquad \rhoh(\ell,m)=0$;
\item $R_5: \qquad \rhoh(\ell,m)\in \bT_1$.
\end{itemize}
Furthermore, we let $R_0=\{(\ell,\ell)\mid \ell \in \cL(q^2)\}$ denote the diagonal relation. In addition, for
$\gre,\gf\in\{-1,1\}$ and $k=0, \ldots, 5$, we let $R_k(\gre,\gf)$ denote the restriction of $R_k$ to
$\cL_\gre(q^2)\times \cL_\gf(q^2)$. For later use, we define sets $\bR_i$ for $i=1,2,\ldots, 5$ by
\[ \bR_1=\bS_0^*, \qquad \bR_2 = \bS_1, \qquad \bR_3 = \bT_0\setminus\bF_q, \qquad
\bR_4=\{0\}, \qquad \bR_5 = \bT_1.\]
Also, we let $r_i=|\bR_i|$ for $i=1,\ldots, 5$, so that
\[ r_1 = (q-2)/2, \qquad r_2 = q/2, \qquad r_3 = q(q-2)/2, \qquad r_4 = 1, \qquad
r_5 = q^2/2.\] Note that for distinct lines $\ell,m\in\cL(q^2)$ and for $k=1,\ldots,5$, we have that $(\ell,m)\in
R_k$ if and only if $\rhoh(\ell,m)\in\bR_k$.

Not all the relations $R_k(\gre,\gf)$ are nonempty.
\begin{lem}\label{Rfus-nonempty}
We have that $R_k(\gre,\gf)$ is nonempty only if\\
(i) $k\in\{0,1,2,3\}$, $\gre=\gf$;\\
(ii) $k=4$, $\gre=\gf=1$;\\
(iii) $k=5$, $\gre\neq \gf$.
\end{lem}
\begin{proof}
Since $\bT_0^*=\bR_1\cup \bR_2\cup \bR_3$, $\bR_4=\{0\}$, and $\bR_5=\bT_1$, this
follows directly from Lemma~\ref{rhoh-type}.
\end{proof}

Each of the relations $R_k(\gre,\gf)$ is a fusion (a union) of relations of the coherent configuration $\cR(q^2)$
on $\cL(q^2)$. We want to show that this fusion, which we will denote by $R(q^2)$, in fact defines a new coherent
configuration. Since this coherent configuration is again weakly symmetric, the restrictions to both
$\cL_{+}(q^2)$ and $\cL_{-}(q^2)$ define fusions $H(q^2)$ and $E(q^2)$ of the hyperbolic and elliptic association
schemes $\cH(q^2)$ and $\cE(q^2)$ defined earlier. In fact, the 4-class fusion $H(q^2)$ of the hyperbolic
association scheme is not new: it was first conjectured to be an association scheme in \cite{CD} and later this
conjecture was proved in \cite{note} (by a direct proof), in \cite{eehx} (by a geometric argument), and in \cite{banstud}
(using characters). The 3-class fusion $E(q^2)$ of the elliptic scheme seems to be new. We will prove these 
fusion results by determining the valencies and the intersection parameters of the fusion.

For all $i,j,k\in\{1,\ldots, 5\}$, for all $c\in\bR_k$, for all $\gre,\gf\in\{-1,1\}$, and for all lines $\ell,m$
with $(\ell,m)\in\cR_c(\gre,\gf)$ (so $\ell\neq m$), we define \beql{vi-def} v_i(\ell)=v_i(\gre)=\#
\{n\in\cL(q^2)\setminus\{\ell\} \mid \rhoh(\ell,n)\in\bR_i\}, \eeql and \beql{pijk-def}
p_{i,j}(\ell,m)=p_{i,j}^k(\gre)= \# \{n\in\cL(q^2)\setminus\{\ell,m\} \mid
        \mbox{$\rhoh(\ell,n)\in\bR_i$ and $\rhoh(n,m)\in\bR_j$}\}.
\eeql
Our aim in this section is to show the following.

\begin{teor}\label{fus-parms}
The numbers $v_i(\gre)$ and $p^k_{i,j}(\gre)$ are well-defined, that is, they do not depend on the particular
choice of the lines $\ell$ and $m$. As a consequence, the relations $R_k(\gre,\gf)$ constitute a coherent
configuration $R(q^2)$ which is a fusion of $\cR(q^2)$. The numbers $v_i(\gre)$ are the valencies of $R(q^2)$;
their values are $v_i(\gre)=r_i(q^2-\gre)$ if $i\neq 4$ and $v_4(\gre)=2(q^2-1)\delta_{\gre,1}$. The numbers
$p^k_{i,j}(\gre)$ are the intersection parameters of $R(q^2)$; their values are given in Tables~\ref{Tpi1} to
\ref{Tpi5}.
\end{teor}

\begin{table}[p]
\tiny
\centering
\begin{tabular}{|c||ccccc|}
\hline
$p^1_{i,j}(\gre)$ & 1 & 2 & 3 & 4 & 5 \\
\hline
1 & $(1+\gre)q^2/2-(4+5\gre)q/2+2+4\gre$ & $q(q-2(\gre+1))/2$   & $(q^2-(4+\gre)q+4(1+\gre))q/2$  & $2(q-3)\delta_{\gre,1}$ &0 \\
2 & $q(q-2(\gre+1))/2$         & $q((1+\gre)q-\gre)/2$   & $q^2(q-(2+\gre))/2$  & $2q\delta_{\gre,1}$     &0 \\
3 & $(q^2-(4+\gre)q+4(1+\gre))q/2$    & $q^2(q-(2+\gre))/2$     & $(q^3-4q^2+(4-\gre)q +2\gre)q/2$ & $2q(q-2)\delta_{\gre,1}$      &0 \\
4 &  $2(q-3)\delta_{\gre,1}$           & $2q\delta_{\gre,1}$      & $2q(q-2)\delta_{\gre,1}$ &$4\delta_{\gre,1}$ &0 \\
5 &         0     &    0        &       0              &     0         & $v_5(\gre)$\\
\hline
\end{tabular}
\caption{\label{Tpi1} Intersection numbers $p^1_{i,j}(\gre)$.}
\end{table}

\begin{table}[p]
\tiny
\centering
\begin{tabular}{|c||ccccc|}
\hline
$p^2_{i,j}(\gre)$ & 1 & 2 & 3 & 4 & 5 \\
\hline
1 & $(q/2-1)(q-2(\gre+1))$  & $(1+\gre)q^2/2-(3\gre+2)q/2+\gre$  & $(q^2-(4+\gre)q+2(\gre+2)q/2$ &$2(q-2)\delta_{\gre,1}$&0 \\
2 & $(1+\gre)q^2/2-(3\gre+2)q/2+\gre$ & $q(q/2-\gre)$            & $(q^2-(2+\gre)q+2\gre)q/2$ &$2(q-1)\delta_{\gre,1}$&0 \\
3 & $(q^2-(4+\gre)q+2(\gre+2)q/2$ &$(q^2-(2+\gre)q+2\gre)q/2$ & $(q^3-4q^2+(4-\gre)q+2\gre)q/2$ & $2q(q-2)\delta_{\gre,1}$ &0 \\
4 & $2(q-2)\delta_{\gre,1}$ & $2(q-1)\delta_{\gre,1}$   & $2q(q-2)\delta_{\gre,1}$ & $4\delta_{\gre,1}$  &0 \\
5 &       0       &    0        &      0               &         0    & $v_5(\gre)$\\
\hline
\end{tabular}
\caption{\label{Tpi2} Intersection numbers $p^2_{i,j}(\gre)$.}
\end{table}

\begin{table}[p]
\tiny
\centering
\begin{tabular}{|c||ccccc|}
\hline
$p^3_{i,j}(\gre)$ & 1 & 2 & 3 & 4 & 5 \\
\hline
1 & $(q^2-(\gre+4)q+4(\gre+1))/2$  & $q(q-\gre-2)/2$  & $(q/2-1)(q^2-2q-\gre)$ & $2(q-2)\delta_{\gre,1}$ &0\\
2 & $q(q-\gre-2)/2$                & $q(q-\gre)/2$    & $(q^2-2q-\gre)q/2$ &$2q\delta_{\gre,1}$      &0\\
3 & $(q/2-1)(q^2-2q-\gre)$ & $(q^2-2q-\gre)q/2$ & $(q^3-4q^2+(4-3\gre)q+8\gre)q/2$ & $2r\delta_{\gre,1}$          &0\\
4 & $2(q-2)\delta_{\gre,1}$ & $2q\delta_{\gre,1}$  & $2r\delta_{\gre,1}$       & $4\delta_{\gre,1}$           &0\\
5 &      0        &   0         &        0             &       0        & $v_5(\gre)$\\
\hline
\end{tabular}
\caption{\label{Tpi3} Intersection numbers $p^3_{i,j}(\gre)$. Here $r=q^2-2q-1$.}
\end{table}

\begin{table}[p]
\footnotesize
\centering
\begin{tabular}{|c||ccccc|}
\hline
$p^4_{i,j}(1)$ & 1 & 2 & 3 & 4 & 5 \\
\hline
1 & $(q-2)(q-3)/2$ & $q(q-2)/2$   &  $ q(q-2)^2/2$        & $q-2$            & 0 \\
2 & $q(q-2)/2 $    & $q(q-1)/2$   &  $ q^2(q-2)/2$        & $q$              & 0 \\
3 & $q(q-2)^2/2$   & $q^2(q-2)/2$ & $q(q-2)r/2$  & $q(q-2)$  & 0 \\
4 & $q-2$          & $q$          &   $q(q-2)$            & $q^2-1$   & 0 \\
5 & 0              &  0           & 0                     & 0 & $v_5(1)$\\
\hline
\end{tabular}
\caption{\label{Tpi4} Intersection numbers $p^4_{i,j}(1)$; here $r=q^2-2q-1$.}
\end{table}

\begin{table}[p]
\footnotesize
\centering
\begin{tabular}{|c||ccccc|}
\hline
$p^5_{i,j}(\gre)$ & 1 & 2 & 3 & 4 & 5 \\
\hline
1 & 0 & 0   & 0   & 0  &  $v_1(\gre)$ \\
2 & 0 & 0 &  0    & 0   & $v_2(\gre)$ \\
3 & 0 & 0 & 0 & 0  & $v_3(\gre)$  \\
4 & 0 & 0 & 0 & 0  & $v_4(\gre)$ \\
5 &$v_1(-\gre)$ &$v_2(-\gre)$ &$v_3(-\gre)$ &$v_4(-\gre)$ &0\\
\hline
\end{tabular}
\caption{\label{Tpi5} Intersection numbers $p^5_{i,j}(\gre)$.}
\end{table}

We will prove this theorem through a sequence of lemmas. First note that we do not need to compute {\em all\/}
the intersection parameters since we have the following:

\begin{lem} \label{sym-sum}
(i) We have that $p^k_{i,0}(\gre)=p^k_{0,i}(\gre)=\delta_{k,i}$ and
$p^0_{i,j}(\gre)= v_i(\gre)\delta_{i,j}$.\\
(ii) If one of $p^k_{i,j}(\gre)$ and $p^k_{j,i}(\gre)$ exists, then so does the other and
$p^k_{j,i}(\gre)=p^k_{i,j}(\gre)$ if $k=1, \ldots, 4$ and
$p^5_{j,i}(\gre)=p^5_{i,j}(-\gre)$.\\
(iii) If four of the numbers $p^k_{i,j}(\gre)$ for $j=1,\ldots, 5$ exist, then so does
the fifth, and
\[ \sum_{j=0}^5 p^k_{i,j}(\gre) = v_i(\gre).\]
\end{lem}
\begin{proof}
Part (i) and (iii) are trivial and part (ii) follows from Lemma~\ref{pi-sym}.
\end{proof}

Our next result justifies all the zero entries in these tables.

\begin{lem}\label{pijk-nz} We have that $p^k_{i,j}(\gre)=0$ for $i,j,k=1, \ldots, 5$
in the following cases.\\
(i) One or three of  $i,j,k$ are equal to 5.\\
(ii) $k=4$ and $\gre=-1$.\\
(iii) $i=4$ or $j=4$, $k=1,2,3$, and $\gre=-1$.
\end{lem}
\begin{proof}
Direct consequence of Lemma~\ref{Rfus-nonempty}.
\end{proof}

As a consequence of Lemmas~\ref{sym-sum} and \ref{pijk-nz} we immediately obtain the intersection parameters
$p^k_{i,j}(\gre)$ for $k=5$.
\begin{teor}\label{int-5}
The numbers $p^5_{i,j}(\gre)$ exist and are as in Table~\ref{Tpi5}.
\end{teor}

To determine the remaining intersection parameters, we proceed as follows. For $c\in\bF_q$, for
$\gre\in\{-1,1\}$, and for $A,B\subseteq \bFqt$, we define
\[ \pi_{A,B}^c(\gre)=\sum_{a\in A, b\in B} \pi^c_{a,b}(\gre). \]
Note that by Lemma~\ref{pi-p}, if the numbers $p^k_{i,j}(\gre)$ exist, then for all $c\in \bR_k$ we have
\beql{int-pi-rel} p^k_{i,j}(\gre)=\sum_{a\in \bR_i, b\in \bR_j} p^c_{a,b}(\gre)=
\pi^c_{\bR_i,\bR_j}(\gre)-\delta_{4,i}\delta_{j,k} - \delta_{4,j}\delta_{i,k}.\eeql So to prove our claim we have
to compute the numbers $\pi^c_{\bR_i,\bR_j}(\gre)$ and show that they do not depend on the choice of $c$ in
$\bR_k$. We need the following simple results.

\begin{lem}\label{pi-zAB}Let $f,g\in\bF_2$ and $\gre\in\{-1,1\}$.
If $A\subseteq \bT_f^*$ and $B\subseteq \bT_g^*$, then
\[\pi^0_{\{0\},\{0\}}(\gre) = \delta_{\gre,1} (q^2+1), \]
\[ \pi^0_{\{0\},B}(\gre) = \pi^0_{B,\{0\}}(\gre) =
              2 |B| \delta_{\gre,1} \delta_{g,0},
\]
and
\[ \pi^0_{A,B}(\gre) =
              (2 |A||B|-|A\cap B| ) \delta_{\gre,1} \delta_{f,g}.
\]
\end{lem}

\begin{proof}
Direct consequence of Theorem~\ref{pi-z}.
\end{proof}

\begin{lem}\label{pi-nzzB}Let $g,h\in\bF_2$ and $\gre\in\{-1,1\}$.
If $B\subseteq \bT_g^*$ and $c\in\bT_h^*$, then
\[\pi^c_{\{0\},\{0\}}(\gre) = 4 \delta_{\gre,1} \delta_{h,0} \]
and
\[ \pi^c_{\{0\},B}(\gre) = \pi^c_{B,\{0\}}(\gre) =
              \delta_{g,h}\biggl((1-2\delta_{\gre,1}) \delta_{c\in B}
                +4 |B| \delta_{\gre,1}\biggr).
\]
\end{lem}
\begin{proof}
Direct consequence of Corollary~\ref{pi-nzz}.
\end{proof}

Using the above results, we now  determine the intersection parameters involving the relation~$R_4$.

\begin{teor}\label{int-k=4}
The intersection parameters $p^4_{i,j}(1)$ exist and are as in Table~$\ref{Tpi4}$.
\end{teor}

\begin{proof}
According to (\ref{int-pi-rel}) and Lemma~\ref{pi-zAB} we have that
\beqa p^4_{4,j}(1) &=& \pi^0_{\{0\},\bR_j}(1) - 2\delta_{j,4}\\
    &=&
    \left\{ \begin{array}{ll}
               2r_j, &\mbox{if $j\in\{1,2,3\}$}; \\
               q^2-1, & \mbox{if $j=4$}.
            \end{array}
    \right.
\eeqa This gives the values for $p^4_{4,j}$ as in Table~\ref{Tpi4}.

Next we consider $p^4_{i,j}(1)$ with $i,j\in\{1,2,3,5\}$, where we assume that either $i=j=5$ or
$i,j\in\{1,2,3\}$. According to (\ref{int-pi-rel}) and Theorem~\ref{pi-z}, we have that
\beqa p^4_{i,j}(1) &=& \pi^0_{\bR_i,\bR_j}(1)\\
    &=&
    \left\{ \begin{array}{ll}
               2r_ir_j, &\mbox{if $i\neq j$}; \\
               (r_i(2(r_i-1)+1)=r_i(2r_i-1), & \mbox{if $i=j$}.
            \end{array}
    \right.
\eeqa
In view of Lemma~\ref{sym-sum}, this is sufficient information to obtain the
remaining values for $p^4_{i,j}(1)$ in Table~\ref{Tpi4}.
\end{proof}

\begin{teor}\label{int-i=4}
The intersection parameters $p^k_{4,j}(1)$ and $p^k_{i,4}(1)$ exist and are as stated in Theorem~$\ref{fus-parms}$.
\end{teor}
\begin{proof}
If $j=k=5$ or if $j,k\in \{1,2,3\}$, and if $c\in \bR_k$,
then using (\ref{int-pi-rel}) and Lemma~\ref{pi-nzzB} we find that
\beqa p^k_{4,j}(1) &=& \pi^c_{\{0\},\bR_j}(1) - \delta_{j,k} \\
    &=&
    \left\{ \begin{array}{ll}
               4r_j, &\mbox{if $j\neq k$}; \\
               4r_j-2, & \mbox{if $j=k$}.
            \end{array}
    \right.
\eeqa
This produces the values for $p^k_{4,j}(1)$ as claimed. 
The other values follow from Lemma~\ref{sym-sum}.
\end{proof}

To complete our determination of the intersection parameters, we compute the numbers $\pi^c_{\bS_i,\bS_j}(\gre)$
for $i,j\in\{0,1\}$ and $c\in\bFqt^*$. Note that both $\bS^*_i$ and $\bS^*_j$ are one of $\bR_1,\bR_2$. Since
\beql{pi-dif} \pi^c_{\bS_i,\bS_j}(\gre) = \pi^c_{\bS^*_i,\bS^*_j}(\gre) + \delta_{i,0}\pi^c_{\{0\},\bS^*_j}(\gre)
+ \delta_{j,0}\pi^c_{\bS^*_i,\{0\}}(\gre) + \delta_{i,0}\delta_{j,0}\pi^c_{\{0\},\{0\}}(\gre) \eeql and since we
know already the numbers $\pi^c_{\{0\},\bR_t}(\gre)$, $\pi^c_{\bR_s,\{0\}}(\gre)$, and
$\pi^c_{\{0\},\{0\}}(\gre)$ with $s,t\in\{1,2\}$, knowing $\pi^c_{\bS_i,\bS_j}(\gre)$ will enable us to find the
numbers $\pi^c_{\bR_s,\bR_t}(\gre)$.

We need some preparation. For $r\in\bF_q$, let us define
\[ \bS_r:= \{ x\in \bFqt \mid \Tr_\bFq(x):= x+x^2 + \cdots + x^{q/2} =r \} \]
and
\[ \bG_r := \{ x\in \bFqt \mid x^q + x = r\}. \]

We will use the following properties.

\begin{lem}\label{SG-con} The above definitions of $\bS_0$ and $\bS_1$ coincide with
the definitions of $\bS_0$ and $\bS_1$ given earlier. Moreover, \\
(i) $\bG_0=\bF_q$, the sets $\bG_r$ with $r\in\bFq$ are precisely the additive cosets of $\bFq$ in $(\bFqt,+)$,
and $\bG_r+\bG_s=\bG_{r+s}$ for all $r,s \in \bFq$. We also have that $\bG_r=r\bG_1$ for $r\in\bF_q^*$.\\
(ii) The map $F: x\mapsto x^2+x$ maps $\bG_r$ two-to-one onto $\bS_r$. In particular, we have $|\bS_r|=q/2$, and
the subsets $\bS_r$ with $r\in\bFq$ partition $\bT_0$. Also, the subsets $\bS_r$ with $r\in\bFq$ are the cosets
of $\bS_0$ in $(\bT_0,+)$;
we have that $\bS_r+\bS_s=\bS_{r+s}$ for all $r,s \in \bFq$.\\
(iii) For each $r\in\bFq$ we have $\bG_{r^2+r} = \bS_r \cup \bS_{r+1}$.
\end{lem}
\begin{proof}
Part (i) follows from the fact that the map $x\mapsto x^q+x$ from $\bF_{q^2}$ to itself is $\bF_2$-linear with
kernel $\bFq$ and image $\bFq$.

Next, if $x\in\bFqt$ satisfies $x^q+x=r$, then
\beqa \Tr_{\bF_q}(x^2+x) &=& (x^2+x) + (x^2+x)^2 + \cdots + (x^2+x)^{q/2} \\
    &=& x+x^q=r,
\eeqa
hence $F$ maps $\bG_r$ to $\bS_r$. Note that the image of $\bFqt$ under $F$ is $\bT_0$;
hence part (ii) now follows from the fact that $F$ is
$\bF_2$-linear with kernel $\bF_2\subseteq \bG_0$.

Finally, if $x^q+x=r$, then $F(x)^q+F(x)=r^2+r$, hence $\bS_r$ (and similarly
$\bS_{r+1}$)
are subsets of $\bG_{r^2+r}$. Since $\bS_r$ and $\bS_{r+1}$ are disjoint and both have
size $q/2$, the result follows.

Remark that $\bG_0=\bFq$, hence $\bS_0$ consists of the elements in $\bFq$ with trace zero, and since
$\bG_0=\bS_0\cup \bS_1$, we have $\bS_1=\bFq\setminus \bS_0$. So the definitions of $\bS_0$ and $\bS_1$ coincide
with the ones given earlier in $\bFq$.
\end{proof}

Now to compute $\pi^c_{\bS_i,\bS_j}(\gre)$ for $\gre=(-1)^e$ with $e\in \bF_2$, for
$i,j\in\{0,1\}$ and for $c\in\bFqt^*$,
we start with the expression
\beqa \pi^c_{\bS_i,\bS_j}(\gre) &=& \sum_{a\in \bS_i}\sum_{b\in \bS_j} \pi^c_{a,b}(\gre)\\
    &=& \sum_{a\in \bS_i}\sum_{b\in \bS_j}
    \biggl( \delta_{a+b+c,0} + 2 \sum_\tau |\bT_e \cap \{ac/\tau^2\}| \biggr),
\eeqa where the sum is over all $\tau \in \bFqt^*$ such that $\tau^2+\tau=a+b+c$. (The last equality was obtained
by using Theorem~\ref{pi-nz}.) Obviously, $\pi^c_{\bS_i,\bS_j}(\gre)$ is nonzero only if $a+b+c\in\bT_0$, hence
only if $c\in\bT_0$. So assume that $c\in \bS_k^*$ with $k\in\bFq$. We will write $r=i+j+k$. Now in the above
expression for $\pi^c_{\bS_i,\bS_j}(\gre)$, we first sum over $b\in  \bS_j$. If $b$ runs through $\bS_j$, then by
Lemma~\ref{SG-con}, part (ii), we have that $a+b+c$ runs through $\bS_r$ and the sum is over all $\tau \in
\bG_r^*$. So we obtain that \beqa \pi^c_{\bS_i,\bS_j}(\gre) &=& \sum_{a\in \bS_i}
    \biggl( \delta_{r,0} + 2 \sum_{\stackrel{\tau^2+\tau \in \bS_r}{\tau\neq0}}
        |\bT_e\cap \{ac/\tau^2\}| \biggr)\\
    &=& \delta_{r,0} q/2  + 2 \sum_{a\in \bS_i}
        \sum_{\tau\in \bG_r^*} |\bT_e\cap \{ac/\tau^2\}| \\
    &=& \delta_{r,0} q/2 + 2 \sum_{\tau\in \bG_r^*}  |\bT_e\cap \bS_i (c/\tau^2)|.
\eeqa
Now we have the following.
\begin{lem}\label{ext-sin} Let $i,e\in\bF_2$ and $\lambda\in\bF_{q^2}^*$. Then
\[
|\bT_e\cap \bS_i\lambda| =
    \left\{ \begin{array}{ll}
              \delta_{e,0} q/2, & \mbox{if $\lambda\in \bG_0^*$}; \\
              \delta_{e,i} q/2, & \mbox{if $\lambda\in \bG_1$}; \\
              q/4, & \mbox{if $\lambda\in \bG_r$ with $r\in\bF_{q^2}\setminus\{0,1\}$}.
            \end{array}
    \right.
\]
\end{lem}
\begin{proof} Since $i\in\bF_2$, we have $\bS_i\subset\bF_q$. For any $s\in \bS_i$, we see that $\lambda s\in \bT_e$
precisely when $e=\Tr_{\bF_{q^2}/\bF_2}(\lambda s)=\Tr_{\bFq}(s(\lambda^q+\lambda))$, that is, when $s\mu \in
\bS_e$, where $\mu=\lambda^q+\lambda$. Now $\mu=0$ if and only if $\lambda\in \bG_0$, in which case for all
$s\in\bS_i$ we have that $s\mu\in \bS_e$ precisely when $e=0$. If $\mu\neq 0$, then the above shows that
$|\bT_e\cap \bS_i\lambda| = |\bS_e\cap \bS_i \mu|$. Note that $\bS_e$ and $\bS_i\mu$ are both hyperplanes of $\bF_q$ 
(considered as a vector space over $\bF_2$), so the size of the intersection equals $\delta_{e,i} q/2$ if $\mu=1$ 
(which occurs precisely when $\lambda \in \bG_1$) and $q/4$ otherwise.
\end{proof}

In order to use this result, given $c\in \bS_k^*$ and $r=k+i+j$ with $i,j\in \bF_2$, we have to determine for how
many $\tau\in \bG_r^*$ we have $c/\tau^2\in \bG_0$, and for how many $\tau\in \bG_r^*$, we have $c/\tau^2\in
\bG_1$. The result is as follows.

\begin{lem}\label{tau-c}
Let $i,j\in\bF_2$, $k \in \bF_q$, $r=k+i+j$, let $\tau\in\bG_r^*$, and let $c\in \bS_k^*$. Write
$c=\gamma^2+\gamma$ with $\gamma\in \bG_k$. Define $\tau_0$, $\tau_1$, and $\tau_2$ by $\tau_0^2=c r/(r+1)$,
$\tau_1^2 = \gamma r$, and $\tau_2^2=(\gamma+1) r$. Then

\noindent (i) $\tau_0$, $\tau_1$, and $\tau_2$ are zero if and only if $r=0$.

\noindent (ii) $\tau_0\in \bG_r$, and for $u=1$ or $2$ we have $\tau_u\in \bG_r$ if and only if $i=j$.

\noindent (iii) $c/\tau^2$ is contained in $\bG_0$ if and only if either $r=0$, or $r\neq0,1$ and $\tau=\tau_0$.

\noindent (iv) $c/\tau^2$ is contained in $\bG_1$ if and only if $r\neq0$ and $\tau\in\{\tau_1,\tau_2\}$.
\end{lem}

\begin{proof}
Obviously, since $c\neq 0$ we have $\gamma\neq 0,1$; hence for $u=0,1,2$, $\tau_u=0$ precisely when $r=0$. Also,
$r\in \bF_q$, and by Lemma~\ref{SG-con} we have $c\in \bS_k\subseteq \bG_{k^2+k}=\bG_{r^2+r}$; hence
\[\tau_0^q+\tau_0=((c^q+c) r/(r+1))^{1/2}= ((r^2+r)r/(r+1))^{1/2} = r,\]
so $\tau_0\in \bG_r$. Also, since $\gamma\in \bG_k$ we have that
\[ \tau_1^q+\tau_1 = (r (\gamma^q+\gamma))^{1/2} =  (rk)^{1/2}, \]
so $\tau_1\in \bG_r$ if and only if $k = r$, that is, $i=j$. The same argument proves the claim for~$\tau_2$.

Finally, let $c/\tau^2 \in \bG_s$ for some $s\in \bF_q$.
We have to determine when $s=0$ and when $s=1$. We saw above that $c\in \bG_{r^2+r}$,
so by definition, we have that
\beqa s &=& (c/\tau^2)^q + c/\tau^2 \\
    &=& (c+r^2+r)/(\tau+r)^2 + c/\tau^2\\
    &=& (cr^2 +(r^2+r) \tau^2)/((\tau (\tau+r))^2.
\eeqa
So firstly, we have $s=0$ if and only if either $r=0$, or $r\neq 0,1$ and
$\tau^2=cr/(r+1)=\tau_0^2$.
Secondly, we can have $s=1$ only if $r\neq 0$. In that case, we have $s=1$ if
\[cr^2+(r^2+r) \tau^2 = (\tau (\tau+r))^2,\]
that is, if $\tau^4+r \tau^2 +c r^2=0$, that is, if $c=(\tau^2/r)^2+\tau^2/r$. So this
happens if $\tau^2/r \in \{\gamma, \gamma+1\}$, that is, if $\tau\in
\{\tau_1,\tau_2\}$.
\end{proof}

\begin{cor}\label{pi-Sij} Let $i,j\in\bF_2$, let $c\in \bS_k^*$ with $k\in \bFq$, and let $e\in\bF_2$. Writing $\epsilon = (-1)^e$, we have that \beq \pi^c_{\bS_i,\bS_j}(\gre) =
    \left\{ \begin{array}{ll}
               q((1+\gre)q-\gre)/2, &\mbox{if $k=i+j$}; \\
               q(2(2\delta_{e,i}-1)+q)/2, & \mbox{if $k=1$ and $i=j$};\\
               q(2(2\delta_{e,i}-1) +q+\gre)/2, & \mbox{if $k\neq0,1$ and $i=j$};\\
               q^2/2, & \mbox{if $k=0$ and $i\neq j$};\\
               q(q+\gre)/2, & \mbox{if $k\neq0,1$ and $i\neq j$}.
            \end{array}
    \right.
\eeq
\end{cor}

\begin{proof}
Let $r=i+j+k$. If we combine Lemma~\ref{ext-sin} and Lemma~\ref{tau-c} with the expression for
$\pi^c_{\bS_i,\bS_j}(\gre)$ just before Lemma~\ref{ext-sin}, we obtain the following. First, if $r=0$, that is,
if $k=i+j$, then
\[ \pi^c_{\bS_i,\bS_j}(\gre) = q/2 + 2 \biggl(
(q-1)\delta_{\gre,1}q/2 \biggr) = q(1+2(q-1)\delta_{\gre, 1})/2.
\]
Next, if $r\neq 0$, that is, if $k\neq i+j$, then $0\not\in \bG_r$, and we obtain that
\beqa
\pi^c_{\bS_i,\bS_j}(\gre) &=& 2 \biggl( \delta_{k\neq 0,1} \delta_{\gre,1} q/2 + 2 \delta_{i=j}\delta_{e,i} q/2 +
(q-\delta_{k\neq 0,1} -2 \delta_{i=j}) q/4
\biggr) \\
    &=&
\biggl( \delta_{k\neq 0,1} \gre  +
2\delta_{i=j}(2\delta_{e,i}-1)  + q
\biggr) q/2,
\eeqa
from which the other expressions follow.
\end{proof}

Now we use (\ref{int-pi-rel}) and (\ref{pi-dif}) to compute the intersection numbers
$p^t_{r,s}(\gre)$ with $r,s \in \{1,2\}$.
For $i,j\in\bF_2$ and $k\in\bFq$, define
\[ \theta(i,j,k)(\gre)=
\delta_{i,0}\pi^c_{\{0\},S^*_j}(\gre)
+ \delta_{j,0}\pi^c_{S^*_i,\{0\}}(\gre)
+ \delta_{i,0}\delta_{j,0}\pi^c_{\{0\},\{0\}}(\gre),
\]
where $c\in \bS_k$ and $c\neq 0$.

\begin{lem} \label{theta}
Write $\gre=2\delta_{e,0}-1=(-1)^e$.
We have that
\[ \theta(i,j,k)(e)=-\gre(\delta_{i,0}\delta_{j,k}+\delta_{j,0}\delta_{i,k}) +
\delta_{e,0} (2q(\delta_{i,0}+\delta_{j,0}) -4\delta_{i,0}\delta_{j,0}). \]
\end{lem}
\begin{proof}
Direct application of Lemma~\ref{pi-nzzB}.
\end{proof}

\begin{teor}\label{prst-nz}
The intersection parameters $p^t_{r,s}(\gre)$ with $r,s \in \{1,2\}$ and $t\in \{1,2,3\}$ exist and are as in
Tables \ref{Tpi1}, \ref{Tpi2}, and \ref{Tpi3}.
\end{teor}

\begin{proof}
Let $r,s \in \{1,2\}$ and $t\in \{1,2,3\}$. Put $i=r-1$ and $j=s-1$ (and consider $i$ and $j$ as elements of
$\bF_2$). If $t=1$, then we take $k=0$; if $t=2$, then we take $k=1$; and if $t=3$, then take $k$ to be any
element in $\bF_q\setminus\{0,1\}$. Finally, let $c\in \bS_k$. Then according to (\ref{int-pi-rel}) and
(\ref{pi-dif}), we have that
\[ p^t_{r,s}(\gre) = \pi^c_{\bS_i,\bS_j}(\gre) - \theta(i,j,k)(e),\]
where $\gre=(-1)^e$. Now we can use Corollary~\ref{pi-Sij} and Lemma~\ref{theta}, firstly to see that the
expression at the right-hand side indeed only depends on $t$ and not on the actual value of $k$ and $c$, and
secondly to compute the value of the intersection parameters $p^t_{r,s}(\gre)$. In this way, we obtain the values
as announced in the theorem.
\end{proof}

Now we can use Lemma~\ref{sym-sum} to find the remaining intersection numbers in
Tables \ref{Tpi1}, \ref{Tpi2} and \ref{Tpi3}.
So we have proved the following.

\begin{teor}\label{prst}
The intersection parameters $p^t_{r,s}(\gre)$ with $t\in\{1,2,3\}$ exist and are as in Tables \ref{Tpi1},
\ref{Tpi2}, and \ref{Tpi3}.
\end{teor}

This completes the proof of Theorem~\ref{fus-parms}.

For the sake of completeness, we mention that the $P$- and $Q$-matrix of the elliptic fusion scheme are given by
\beql{Pell} P=
\left(
\begin{array}{cccc}
1   & (q-2)(q^2+1)/2  & q(q^2+1)/2  & q(q-2)(q^2+1)/2  \\
1   & -(q-1)(q-2)/2   & -q(q-1)/2   & q(q-2)           \\
1   & -(q^2-q+2)/2    & q(q+1)/2    & -q               \\
1   & q-1             & 0           & -q
\end{array}\right)
\eeql
and \beql{Qell} Q=
\left(
\begin{array}{cccc}
1   & q(q^2+1)/2  & (q-2)(q^2+1)/2  & q(q-2)(q^2+1)/2  \\
1   & -q(q-1)/2   & -(q^2-q+2)/2   &  q(q-1)           \\
1   & -q(q-1)/2    & (q-2)(q+1)/2    & 0               \\
1   & q           & -1           & -q
\end{array}\right).
\eeql
The $P$- and $Q$-matrix of the hyperbolic fusion scheme can be found in \cite{eehx}.
\section{\label{fus-plus}Further fusions}
From the values of the intersection parameters $p^k_{i,j}(\gre)$ as computed in the previous section we
immediately see that a further fusion of the relations $R_1(\gre,\gre)$ and $R_2(\gre,\gre)$ for $\gre=1$ and
$\gre=-1$ produces another weakly symmetric coherent configuration, and thus a further 3-class association scheme
on the hyperbolic lines and 2-class association scheme (that is, a strongly regular graph) on the elliptic lines;
the intersection parameters are given in Tables \ref{TFpi12}, \ref{TFpi3}, \ref{TF4ij}, and \ref{TF5ij} below.

\begin{table}[p]
\tiny
\centering
\begin{tabular}{|c||cccc|}
\hline
$p^{\{1,2\}}_{r,s}(\gre)$ & \{1,2\} & 3 & 4 & 5 \\
\hline
\{1,2\} & $(2+\gre)q^2-(4+5\gre)q+2+4\gre$ & $q(q^2-(\gre+3)q +2(1+\gre))$ & $2(2q-3)\delta_{\gre,1}$ &0 \\
3 & $q(q^2-(\gre+3)q +2(1+\gre))$ & $(q^3-4q^2+(4-\gre)q+2\gre)q/2$ & $2q(q-2)\delta_{\gre,1}$ &0 \\
4 &  $2(2q-3)\delta_{\gre,1}$      & $2q(q-2)\delta_{\gre,1}$ &$4\delta_{\gre,1} $ &0 \\
5 &         0       &       0              &     0           & $v_5(\gre)$\\
\hline
\end{tabular}
\caption{\label{TFpi12} Intersection numbers $p^{\{1,2\}}_{r,s}(\gre)$.}
\end{table}

\begin{table}[p]
\tiny
\centering
\begin{tabular}{|c||cccc|}
\hline
$p^3_{r,s}(e)$ & \{1,2\} & 3 & 4 & 5 \\
\hline
\{1,2\} & $2q^2-(2\gre+4)q+2(\gre+1)$  & $q^3-3q^2-(\gre-2)q+\gre$  & $4(q-1)\delta_{\gre,1}$ &0\\
3 & $q^3-3q^2-(\gre-2)q+\gre$ & $(q^3-4q^2+(4-3\gre)q+8\gre)q/2$ & $2r\delta_{\gre,1}$  &0\\
4 & $4(q-1)\delta_{\gre,1}$ &  $2r\delta_{\gre,1}$       & $4\delta_{\gre,1}$           & 0\\
5 &      0          &        0             &       0           & $v_5(\gre)$\\
\hline
\end{tabular}
\caption{\label{TFpi3} Intersection numbers $p^3_{r,s}(\gre)$. Here $r=q^2-2q-1$.}
\end{table}

\begin{table}[p]
\footnotesize
\centering
\begin{tabular}{|c||cccc|}
\hline
$p^4_{i,j}(1)$ & \{1,2\} & 3 & 4 & 5 \\
\hline
\{1,2\} & $2q^2-5q+3$ & $ q(q-1)(q-2)$  & $2(q-1)$   & 0 \\
3 & $q(q-1)(q-2)$   &  $q(q-2)r/2$  & $q(q-2)$  & 0 \\
4 & $2(q-1)$          &   $q(q-2)$            & $q^2-1$   & 0 \\
5 & 0              &  0           & 0                     & $v_5(1)$\\
\hline
\end{tabular}
\caption{\label{TF4ij} Intersection numbers $p^4_{i,j}(1)$; here $r=q^2-2q-1$.}
\end{table}

\begin{table}[p]
\footnotesize
\centering
\begin{tabular}{|c||cccc|}
\hline
$p^5_{i,j}(\gre)$ & \{1,2\} &  3 & 4 & 5 \\
\hline
\{1,2\} &  0   & 0   & 0  &  $v_1(\gre)+v_2(\gre)$ \\
3 & 0 & 0  & 0  & $v_3(\gre)$  \\
4 & 0 & 0  & 0  & $v_4(\gre)$ \\
5 &$v_1(-\gre)+v_2(-\gre)$ &$v_3(-\gre)$ &$v_4(-\gre)$ &0\\
\hline
\end{tabular}
\caption{\label{TF5ij} Intersection numbers $p^5_{i,j}(\gre)$.}
\end{table}

Finally, we see that a further fusion of $R_1(1,1)\cup R_2(1,1)$ with $R_4(1,1)$ again produces a weakly
symmetric coherent configuration, and thus a 2-class association scheme (that is, a strongly regular graph) on
the hyperbolic lines. Some of the intersection parameters of this further fusion are given in Tables
\ref{TFFpi12} and \ref{TFFpi3}.

\begin{table}[p]
\footnotesize
\centering
\begin{tabular}{|c||ccc|}
\hline
$p^4_{i,j}(1)$ & \{1,2,4\} & 3 & 5 \\
\hline
\{1,2,4\} & $3q^2-q-2$ & $ q^2(q-2)$  &  0 \\
3 & $q^2(q-2)$   &  $q(q-2)r/2$  &  0 \\
5 & 0              &  0        & $v_5(1)$\\
\hline
\end{tabular}
\caption{\label{TFFpi12} Intersection numbers $p^{\{1,2,4\}}_{i,j}(1)$;
here $r=q^2-2q-1$.}
\end{table}

\begin{table}[p]
\tiny
\centering
\begin{tabular}{|c||ccc|}
\hline
$p^3_{r,s}(e)$ & \{1,2,4\} & 3 & 5 \\
\hline
\{1,2,4\} & $2q(q+2)$  & $(q+1)r$  & 0\\
3 & $(q+1)r$ & $(q^3-4q^2+q+8)q/2$ & 0\\
5 &      0   &   0               & $v_5(1)$\\
\hline
\end{tabular}
\caption{\label{TFFpi3} Intersection numbers $p^3_{r,s}(1)$. Here $r=q^2-2q-1$.}
\end{table}

So in this way we obtain two strongly regular graphs, with parameters $(v,k,\lambda,\mu)$ where
$v=q^2(q^2+\gre)/2$, $k=(q^2-\gre)(q+\gre)$, $\lambda=2(q^2-1) +\gre q(q-1)$, and $\mu=2q(q+\gre)$. Graphs with
these parameters were first described by R.~Metz for $\gre=-1$ (the elliptic case) and by Brouwer and Wilbrink
for $\gre=1$ (the hyperbolic case), see \cite[Section~7]{bl}. The two constructions were further generalized by
Wilbrink. For $q=4$, the ``elliptic'' graph was obtained and conjectured to be a Metz graph in \cite[p.~83]{henkthesis}.

In a forthcoming paper \cite{hh-iso} it will be shown among other things that the strongly regular graphs
obtained above are in fact {\em isomorphic\/} to the Metz graphs (for $\gre=-1$) or the Brouwer-Wilbrink graphs
(for $\gre=1$). For $\gre=1$ this has already been proved in \cite{eehx}; for $\gre=-1$ this was first
conjectured in \cite{henkthesis} for $q=4$, and proved for the first time in \cite{hh-iso}.

\vskip3pt \noindent{\bf Acknowledgements:} The second author thanks Philips Research Eindhoven, the Netherlands,
where part of this work was carried out. The research of the second author is supported in part by NSF grant DMS
0400411. Also, the authors thank Ludo Tolhuizen and Sebastian Egner for some useful comments.

\end{document}